\documentclass[reqno,12pt]{amsart}
\usepackage{amsmath,amsthm,amssymb,amsfonts,amscd}
\usepackage{xypic}
\usepackage{xcolor}
\theoremstyle{plain} 
	\newtheorem{thm}{Theorem}[section]
	\newtheorem*{thm*}{Theorem}
	\newtheorem{cor}[thm]{Corollary}
	\newtheorem{lem}[thm]{Lemma}
	
	\newtheorem{prop}[thm]{Proposition}
	
	\newtheorem*{conj*}{Conjecture}
	
\theoremstyle{definition}
	\newtheorem{defn}[thm]{\rm{Definition}}%[section]

\theoremstyle{remark}
	\newtheorem{rmk}[thm]{\rm{Remark}}
	
	\newtheorem*{pf}{\rm{Proof}}
\numberwithin{equation}{section}
\def\CC{{\mathbb C}}

\def\RR{{\mathbb R}}
\def\ZZ{{\mathbb Z}}

\def\p{\partial }

\def\multi{{\mathbb{Z}^n_{\geq 0}}}
\def\one{{1 \leq \alpha \leq n}}
\def\two{{1 \leq \alpha,\,\beta \leq n}}
\def\Der{{Der(S^W)}}
\def\rest{{_{L^{\perp}}}}

\newcommand{\bp}{\begin{pmatrix}}
\newcommand{\ep}{\end{pmatrix}}

\numberwithin{equation}{section}

\newcounter{CounterEQUlabel}
\newcommand{\EQUlabel}[1]{\label{#1}
	\ifcase \theCounterEQUlabel
		\relax
	\or
		\hspace{1em}\mbox{\tiny$\langle$\rmfamily#1$\rangle$}
		\index{zzz#1@#1}
	\fi }	
	\newcounter{CounterEQUref}
	\newcounter{CounterEQUpageref}
	\newcommand{\EQUref}[1]{
		\ifcase \theCounterEQUref     \relax   \or {\tiny[#1]}\,\fi
		\ifcase \theCounterEQUpageref (\ref{#1}) \or (\ref{#1}\,(p.\pageref{#1})) \fi}

\begin{document}
\title{Good basic invariants for elliptic Weyl groups and Frobenius structures}
\footnote{2010 Mathematics Subject Classification. Primary 32G20; Secondary 32N15.}
\date{\today}
\author{Ikuo Satake}
\address{Faculty of Education, Bunkyo University, 
3337 Minamiogishima Koshigaya, Saitama, 343-8511, Japan}
\email{satakeikuo@gmail.com}
\begin{abstract}
In this paper, we define a set of good basic invariants for 
the elliptic Weyl group for the elliptic root system. 
For an elliptic root system of codimension $1$, 
we show that a set of good basic invariants 
gives a set of flat invariants obtained by Saito 
and that Taylor coefficients of the good basic invariants 
give the structure constants of the multiplication of the 
Frobenius structure obtained by the author. 
\end{abstract}
\maketitle
%%%%%%%%%%%%%%%%%%%%%%%%%%%%%%%%%%%%%%%%%%%%%%%%%%%%%%%%%%%%%%%%
\section{Introduction}

\subsection{Aim and results of the paper}
Let $W$ be an elliptic Weyl group 
defined as a reflection group 
for an elliptic root system. 
Let $Y$, $H$ be domains 
where $H$ is isomorphic to the upper half plane 
and $\pi:Y \to H$ is an affine bundle whose fiber 
is isomorphic to $\CC^{l+1}$. 
The group $W$ acts on each fiber of $\pi$. 

Let $\widetilde{c} \in W$ be a hyperbolic Coxeter transformation 
defined in \cite{extendedI}. 
It is known that $\widetilde{c}$ is not semi-simple and has no fixed points 
on $Y$. 

We take the Jordan decomposition:
$$
\widetilde{c}=\widetilde{c}^{ss}\cdot \widetilde{c}^{unip}, 
$$
where $\widetilde{c}^{ss}$ is a semi-simple part and 
$\widetilde{c}^{unip}$ is a unipotent part. 
We show in Section 7 that we could take
a suitable section $S \subset Y$ 
(which gives isomorphism $S \simeq H$ by the composite morphism 
$S \subset Y \to H$) such that 
\begin{enumerate}
\item every point in $S$ is fixed by the action of $\widetilde{c}^{ss}$,
\item no points of $S$ are contained in the reflection hyperplanes of $W$. 
\end{enumerate}
Then for the $W$-invariants on $Y$, 
we define Taylor expansions along $S$ and 
by using these Taylor expansions, we define a set of good basic invariants 
which is analogous to the cases for the Coxeter groups (see \cite{handai3}). 

In this paper, we define the notion of an admissible triplet 
$(g,\zeta,L)$ (cf. Definition \ref{200310.1}) 
which has the same role as $\widetilde{c}^{ss}(=g)$ and $S(=L^{\perp})$. 
Then we define a set of good basic invariants. 

For the elliptic root systems of codimension $1$, 
we show that a set of good basic invariants 
gives a set of flat invariants obtained by Saito \cite{extendedII} 
and the Taylor coefficients of the good basic invariants 
give the structure constants of the multiplication of the 
Frobenius structure obtained by the author \cite{handai}. 

In the study of the Frobenius structure for the elliptic Weyl groups, 
a characterization (Theorem \ref{300.001}\hskip0.2mm(iii)) 
of the unit field is important. 
We find another characterization (Proposition \ref{200323.3}) of 
the unit field by the space $S(=L^{\perp})$. 
This enables us to find the notion of the 
good basic invariants for elliptic Weyl groups 
and the ones for finite complex reflection groups. 

Here is a brief account of the contents of the paper. 
In Section 2 we remind notions of elliptic 
root systems, elliptic Weyl groups and 
their invariant theory. 
In Section 3 we define an admissible triplet for 
the elliptic root system. 
In Section 4 we define good basic invariants. 
In Section 5 we give properties of Taylor coefficients of 
good basic invariants.  
In Section 6 we construct an admissible triplet. 
In Section 7 we show the uniqueness of good basic invariants 
under suitable assumptions. 
In Section 8 We treat the elliptic root system of type ``codimension one". 
We give a description of the bilinear form in terms of 
the good basic invariants (Theorem \ref{230403.1}). 
In Section 9 we show that the good invariants give a nice description 
of the Frobenius structure which is defined by 
Saito and Satake. 

\subsection{Acknowledgements}
This work is supported in part by 
Grant-in Aid for Challenging Research (Exploratory) 17K18781 
, 
Grant-in-Aid for Scientific Research(C) 18K03281 
and 
Grant-in-Aid for Scientific Research(C) 22K03295 
\vskip1cm

%\newpage

\section{$W$-invariants for elliptic root systems}
We recall the elliptic root systems (cf. Saito \cite{extendedII}). 
We use notations in \cite{handai2} in order to fit our notations to the ones 
of Kac \cite{Kac}. 

\subsection{Elliptic root systems}\label{240311.2}

In this subsection, we define an elliptic root system (cf. \cite{extendedII}). 

Let $l$ be a positive integer. 
Let $F$ be a real vector space of rank $l+2$ with a positive semi-definite 
symmetric bilinear form 
$I:F \times F \to \RR$, whose radical 
$\mathrm{rad}\,I:=\{x \in F\,|\,I(x,y)=0,\forall y \in F\}$ 
is a vector space of rank $2$. 
We put 
$O(F,\mathrm{rad}\,I):=\{g \in GL(F)\,|\,I(gx,gy)=I(x,y)\,\forall x,y \in F,\ 
g|_{\mathrm{rad}\,I}=id.\,\}$. 
For a non-isotropic vector $\alpha \in F$ 
(i.e. $I(\alpha,\alpha)\neq 0$), 
we put $\alpha^{\vee}:=2\alpha/I(\alpha,\alpha)\in F$. 
The reflection $w_{\alpha} \in O(F,\mathrm{rad}\,I)$ with respect to $\alpha$ is 
defined by 
$
w_{\alpha}(u):=u-I(u,\alpha^{\vee})\alpha
\quad
(\forall u \in F)
$. 
\begin{defn}
A set $R$ of non-isotropic elements 
of $F$ is an elliptic root system belonging to $(F,I)$ if it satisfies 
the axioms (i)--(iv). 
\begin{enumerate}
\item The additive group generated by $R$ in $F$, 
denoted by $Q(R)$, 
is a full sub-lattice of $F$.
\item $I(\alpha,\beta^{\vee}) \in \ZZ$ for $\alpha,\beta \in R$. 
\item $w_{\alpha}(R)=R$ for $\forall \alpha \in R$. 
\item If $R=R_1 \cup R_2$, with $R_1 \perp R_2$, then either 
$R_1$ or $R_2$ is void. 
\end{enumerate}
\end{defn}
We have 
$Q(R)\cap\hskip0.2mm\mathrm{rad}\,I \simeq\ZZ^2$. 
We call a 1-dimensional vector space 
$G \subset\hskip0.2mm \mathrm{rad}\,I$ 
satisfying $G\cap Q(R) \simeq \ZZ$, a marking. 
We fix $a,\delta \in F$ s.t. 
$G \cap Q(R)=\ZZ a$ and 
$Q(R) \cap\hskip0.2mm\mathrm{rad}\,I=\ZZ a \oplus \ZZ \delta$. 
Let $(I_R:I) \in \RR_{> 0}$ be the smallest number such that 
$(I_R:I)I$ defines an even lattice structure on $Q(R)$. 
The bilinear form $(I_R:I)I$ is denoted by $I_R$. 

The isomorphism classes of the elliptic root systems 
with markings are classified in \cite{extendedI}. 

\subsection{Hyperbolic extension}

We introduce a hyperbolic extension $(\widetilde{F},\widetilde{I})$ of $(F,I)$, 
i.e. $\widetilde{F}$ is a $(l+3)$-dimensional $\RR$-vector space of 
which contains $F$ as a subspace 
and $\widetilde{I}$ is a symmetric $\RR$-bilinear form 
on $\widetilde{F}$ which satisfies $\widetilde{I}|_F=I$ and 
$\mathrm{rad}\,\widetilde{I}=\RR a$. 
It is unique up to isomorphism. 
We put $O(\widetilde{F},F,\mathrm{rad}\,I):=\{g \in GL(\widetilde{F})\,|\,
\widetilde{I}(gx,gy)=\widetilde{I}(x,y)\,\forall x,y \in \widetilde{F},\,
g(F) \subset F,\,g|_F \in O(F,\mathrm{rad}\,I)\,\}$. 
The natural homomorphism 
$
O(\widetilde{F},F,\mathrm{rad}\,I)\to O(F,\mathrm{rad}\,I)$ is surjective 
and its kernel $K_{\RR}$ is isomorphic to the additive group $\RR$:
\begin{equation}\label{200322.4}
0 \to K_{\RR} \to 
O(\widetilde{F},F,\mathrm{rad}\,I)\to O(F,\mathrm{rad}\,I)
\to 1.
\end{equation}

We fix some notations. 
We take $\Lambda_0 \in \widetilde{F}$ which satisfies 
$\widetilde{I}(\Lambda_0,\delta)=1$ and $\widetilde{I}(\Lambda_0,\Lambda_0)=0$. 
Then we have a decomposition $\widetilde{F}=F\oplus \RR\Lambda_0$. 
\subsection{Elliptic Weyl group}
We define an elliptic Weyl group. 
For $\alpha \in R$, we define a reflection $\widetilde{w}_{\alpha} \in 
O(\widetilde{F},F,\mathrm{rad}\,I)$ by 
$\widetilde{w}_{\alpha}(u)=u-\widetilde{I}(u,\alpha^{\vee})\alpha$ for $u \in \widetilde{F}$.  
We define an elliptic Weyl group $W$ as a group 
generated by $\widetilde{w}_{\alpha}\ (\alpha \in R)$. 
A subgroup $K_{\ZZ}:=W \cap K_{\RR}$ is isomorphic to $\ZZ$. 
\subsection{Domains and Euler field}\label{240.001}

We define two domains:
\begin{eqnarray}
&&Y:=\{x \in \mathrm{Hom}_{\RR}(\widetilde{F},\CC)\,|\,\langle a,x\rangle=-2\pi\sqrt{-1},\ 
		\mathrm{Re}\langle \delta,x\rangle>0\},\label{230330.1}\\
&&H:=\{x \in \mathrm{Hom}_{\RR}(\mathrm{rad}\,I,\CC)\,|\,\langle a,x\rangle=-2\pi\sqrt{-1},\ 
		\mathrm{Re}\langle \delta,x\rangle>0\}.\label{230330.2}
\end{eqnarray}
We have a natural morphism 
\begin{equation}\label{230330.3}
\pi:Y\to H. 
\end{equation}
We remark that 
$H$ is isomorphic to $\mathbb{H}:=\{z\in \CC\,|\,\mathrm{Im}z>0\}$ 
by the function 
$\delta/(-2\pi\sqrt{-1})=\delta/a:H \to \mathbb{H}$. 
We define the left action of $g \in W$ on $Y$ by 
$\langle g\cdot x,\gamma \rangle=\langle x,g^{-1} \cdot \gamma \rangle$ 
for $x \in Y$ and $\gamma \in \widetilde{F}$. 
We remark that the domain $Y$ and the action of $W$ on $Y$ 
are naturally identified with the ones of Kac \cite[p.225]{Kac}. 

Let $F(Y)$ be the space of holomorphic functions on $Y$ 
and $\Omega(Y)$ be the space of holomorphic 1-forms on $Y$. 
We denote by $m_0 \leq \cdots\leq m_l$ the exponents of the 
elliptic root system with marking (see \cite[p22]{extendedII}) 
and we also denote by $m_{max}$ the maximum of the exponents. 

We put $n=l+1$. 
Let $d_n$ be the smallest common denominator for the rational numbers 
$m_i/m_{max}$ ($i=0,\cdots,l$). 
We define the normalized exponents by 
\begin{equation}\label{240311.1}
d_{\alpha}:=m_{\alpha-1}\frac{d_n}{m_{max}}\quad
(\alpha=1,\cdots,n-1).
\end{equation}
We remark that the equation (\ref{240311.1}) holds for the case 
$\alpha=n$ because $m_l=m_{max}$. 
We also remark that $d_n$ is $l_{max}+1$ in \cite[p23]{extendedII}. 

There exists a unique vector field $E$  
\begin{equation}\label{200322.1}
E:\Omega(Y) \to F(Y)
\end{equation}
such that $E(f)=0$ for any $f \in F$ and 
\begin{equation}
E(\Lambda_0)=\frac{(I_R:I)d_n}{m_{max}}. 
\end{equation}
We call $E$ {\it the Euler field}.
%

%%%%%%%%%%%%%%%%%%%%%%%%%%%%%%%%%%%%%%%%%%%%%%%%%%%
\subsection{The algebra of the invariants for the elliptic Weyl group}
In this subsection, we introduce the algebra of the invariants 
for the elliptic Weyl group. 

We put $F({H}):=\{f:{H} \to \CC:holomorphic\}$. 
For $m \in\ZZ$, we put 
\begin{equation}\label{230817.1}
F_m(Y):=\{f\in F(Y)\,|\,Ef=mf\}. 
\end{equation}
The morphism $\pi:Y \to H$ induces $\pi^*:F({H}) \to F_0(Y)$, 
thus $F_0(Y)$-module $F_m(Y)$ is an $F({H})$-module. 

For $m \in \ZZ$ , we put 
\begin{eqnarray}
&&S_m^W:=\{f\in F_m(Y)\,|\,f(g\cdot z)=f(z),\ \forall g \in W\}, \\
&&S^W:=\bigoplus_{m \in \ZZ_{\geq 0}}S^W_{m}. \label{200309.1}
\end{eqnarray}

$S^W$ is a graded $F({H})$-algebra. 

\begin{thm}[\cite{Chevalley3, Chevalley4, 
Chevalley6, Chevalley7, Chevalley8, 
Chevalley2, Chevalley1, Chevalley5}]
The $F({H})$-algebra $S^W$ is generated by a set of algebraically independent homogeneous generators $x^1,\cdots,x^{n}$ $(n=l+1)$ with degrees 
$0 <d_1 \leq d_2 \cdots \leq d_n$ which we call a set of basic invariants. 
\end{thm}
We put 
\begin{equation}
d:=(d_1,\cdots,d_n).
\end{equation}

Then the degree $m$ part $S_m^W$ could be also written as 
\begin{equation}
S_m^W=\{\sum_{b \in \multi} A_b x^b\in S^W\,|\, A_b \in F(H),\ 
d\cdot b=m\},
\end{equation}
where we denote 
\begin{equation}
x^b=(x^1)^{b_1}\cdots (x^n)^{b_n},\quad
d\cdot b=d_1 b_1+\cdots+d_nb_n. 
\end{equation}

%%%%%%%%%%%%%%%%%%%%%%%%%%%%%%%%%%%%%%%%%%%%%%%%%%%
\subsection{Decomposition of the space $Y$}\label{230813.4}

Let $\widetilde{X}$ be a space of 
complementary subspaces of $\mathrm{rad}\,I$ in a vector space $\widetilde{F}$:
\begin{equation}\label{230812.5}
\widetilde{X}:=\{V \subset \widetilde{F}\,|\,
\widetilde{F}=V \oplus \mathrm{rad}\,I\}.
\end{equation}
The space $\widetilde{X}$ is an affine space over 
$\mathrm{Hom}_{\RR}(\widetilde{F}/\mathrm{rad}\,I,\mathrm{rad}\,I)$. 

For the space $Y$ defined in (\ref{230330.1}), 
we define a mapping:
\begin{equation}
f_1:Y \to \widetilde{X}
\end{equation}
by $f_1(y)=\mathrm{ker}\ y$ where we see $y \in Y$ 
as a morphism $y:\widetilde{F} \to \CC$. 
We see that $f_1$ is $O(\widetilde{F},F,\mathrm{rad}\,I)$-equivariant. 
By $\pi:Y \to H$ defined in (\ref{230330.3}), 
we have a mapping:
\begin{equation}\label{230812.1}
(\pi,f_1):Y \to H \times \widetilde{X}
\end{equation}
which is an isomorphism as a real manifold. 

The mapping $f_1$ is not a holomorphic mapping, 
but for any $V \in \widetilde{X}$, the subset $f^{-1}(V)$ 
has a description 
\begin{equation}\label{240314.1}
f_1^{-1}(V)=\{x \in Y\,|\,\langle v,x \rangle=0\ \forall v \in V\}
\end{equation}
which gives the structure of a complex submanifold of $Y$. 
We remark that $f_1^{-1}(V)$ is isomorphic to $H$. 
Then the mapping $f_1$ gives a decomposition of 
$Y$ into complex submanifolds
\begin{equation}\label{230812.4}
Y=\bigsqcup_{V \in \widetilde{X}}f_1^{-1}(V). 
\end{equation}

%If we fix one splitting subspace $V_0 \in \widetilde{X}$ and 
%an isomorphism $V_0 \simeq \RR^n$, 
%we have 
%\begin{subequations}
%\begin{eqnarray}
%\widetilde{X}
%&\simeq& 
%\{x \in \mathrm{Hom}_{\RR}(\widetilde{F},\mathrm{rad}\,I)\,|\,
%x|_{\mathrm{rad}\,I}=id.\} \label{230813.1}\\
%&\simeq& 
%\mathrm{Hom}_{\RR}(V_0,\mathrm{rad}\,I)\label{230813.2}\\
%&\simeq&
%(\mathrm{rad}\,I)^n, \label{230813.3}
%\end{eqnarray}
%\end{subequations}
%where $(\mathrm{rad}\,I)^n$ is a direct sum of $n$-copies of $\mathrm{rad}\,I$. 
%We also have 
%\begin{eqnarray*}
%Y
%&\simeq&
%H \times \mathrm{Hom}_{\RR}(V_0,\CC)\\
%&\simeq&
%H \times \CC^n.
%\end{eqnarray*}
%Then the inverse mapping of $(\ref{230812.1})$ is described as 
%$$
%H \times (\mathrm{rad}\,I)^n
%\to 
%H \times \CC^n
%$$
%for $(x_1,p_1a+q_1\delta,\cdots,p_n a+q_n \delta) \mapsto 
%(x_1,
%p_1 \langle a,x_1 \rangle+q_1\langle \delta,x_1\rangle,
%\cdots,
%p_n \langle a,x_1 \rangle+q_n\langle \delta,x_1\rangle)
%$, where $x_1 \in H$ and $p_i,q_i \in \RR\ (1 \leq i \leq n)$.
%
%If $x \in \widetilde{X}$ corresponds to 
%$(p_1a+q_1\delta,\cdots,p_n a+q_n \delta) \in (\mathrm{rad}\,I)^n$ 
%by $(\ref{230813.1})$, $(\ref{230813.2})$ and $(\ref{230813.3})$, 
%then $f_1^{-1}(x)$ in $(\ref{230812.4})$ is isomorphic to  
%$$
%\{(x_1,z_1,\cdots,z_n) \in H \times \CC^n\,|\,
%z_i=p_i \langle a,x_1 \rangle+q_i\langle \delta,x_1\rangle 
%(1 \leq i \leq n)\}
%$$
%which is a $1$-dimensional complex submanifold, 

%\newpage

%%%%%%%%%%%%%%%%%%%%%%%%%%%%%%%%%%%%%%%%%%%%%%%%%%%%%%%%%%%%%%%%%%%%%%%%%%
\section{Graded algebra}
For an elliptic root system, 
we define a notion of an admissible triplet. 
\subsection{Admissible triplet}
\begin{defn}\label{200310.1}
For $g \in O(\widetilde{F},F,\mathrm{rad}\,I),\ 
\zeta \in \mathbb{C}^*$ and $L \subset \widetilde{F}$, 
we call a triplet 
$(g,\zeta,L)$ {\it admissible} if it 
satisfies the following conditions. 
\begin{enumerate}
\item $g$ is semi-simple and 
\begin{equation}\label{230327.2}
(g-id.)(\widetilde{F}) \subset F. 
\end{equation}
\item $\zeta$ is a primitive $d_n$-th root of unity and 
\begin{equation}\label{200322.7}
g\cdot x^{\alpha}=\zeta^{d_{\alpha}}x^{\alpha}\quad (1 \leq \alpha \leq n)
\end{equation}
for a set of basic invariants $x^1,\cdots,x^n$ with degrees 
$d_1,\cdots,d_n$. 
\item $L$ is a splitting subspace of $\mathrm{rad}\,I$, which is 
$g$-stable and has no roots:
\begin{eqnarray}
&&\widetilde{F}=L \oplus \mathrm{rad}\,I, \label{200323.7}\\
&& g(L) =L, \label{230329.4}\\
&&L \cap R =\emptyset. \label{200322.8}
\end{eqnarray}
We remark that $\mathrm{dim}L=n(=l+1)$. 
\end{enumerate}
\end{defn}
From now on we fix an admissible triplet $(g,\zeta,L)$. 

By Definition \ref{200310.1}\hskip0.2mm(i), the action of $g$ on $L$ is 
also semi-simple. 
We take 
a $\CC$-basis $z^1,\cdots,z^n$ of $L \otimes_{\RR}\CC$ such that 
\begin{equation}\label{230327.1}
g \cdot z^{\alpha}=c_{\alpha}z^{\alpha}\quad (1 \leq \alpha \leq n),
\end{equation}
where $c_1,\cdots,c_n$ are eigenvalues of $g$ on $L$. 
We consider $z^1,\cdots,z^n$ as functions on $Y$ and 
we define $z^0 \in F(Y)$ by 
\begin{equation}\label{200323.12}
z^0(x):=\frac{\langle \delta,x\rangle}{-2\pi\sqrt{-1}} \quad (x \in Y). 
\end{equation}
Then the set $z^0,z^1,\cdots,z^n$ gives a coordinate system of $Y$. 

\begin{defn}
We put
\begin{equation}
L^{\perp}:=\{x \in Y\,|\,\langle l,x\rangle=0\ \forall\,l \in L\,\}.
\end{equation}
\end{defn}
We remark that $L \in \widetilde{X}$ and 
$L^{\perp}=f_1^{-1}(L)$ (see $(\ref{240314.1})$. 
The morphism $L^{\perp} \to H$ induced by $\pi:Y \to H$ is an isomorphism. 
Then we identify the restriction $f|_{L^{\perp}}$ 
of a function $f \in F(Y)$ to $L^{\perp}$ with the function on $H$. 

\begin{prop}\label{230402.1}
For any $q \in L^{\perp}$ and $g$ of the admissible triplet $(g,\zeta,L)$,
 we have 
\begin{equation}
g\cdot q=q.
\end{equation}
\end{prop}
\begin{pf}
$g$ acts on $L$, then $g \cdot q \in L^{\perp}$. 
$\langle g \cdot q,a \rangle 
=\langle q,a \rangle$, 
$\langle g \cdot q,\delta \rangle 
=\langle q,\delta \rangle$, 
then $g \cdot q=q$. 
\qed\end{pf}

We put 
\begin{eqnarray}
&&H_{\alpha}:=\{x \in Y\,|\,\langle \alpha,x\rangle=0\,\}
\end{eqnarray}
for $\alpha \in R$. 
\begin{prop}\label{200322.12}
\begin{enumerate}
\item The space $L^{\perp}$ is regular, i.e. 
\begin{equation}
L^{\perp} \cap \left(
\bigcup_{\alpha \in R} H_{\alpha}
\right)
=\emptyset.
\end{equation}
\item
On the space $L^{\perp}$, the Jacobian matrix 
\begin{equation}\label{200318.10}
\left(
\left.
\frac{\p x^{\alpha}}{\p z^{\beta}}
\right|_{L^{\perp}}
\right)
_{1 \leq \alpha,\beta \leq n}
\end{equation}
is invertible. 
\item The eigenvalues of $g$ on $L \otimes_{\RR} \CC$ are 
$\zeta^{d_{\alpha}}$ $(1 \leq \alpha \leq n)$. 
\end{enumerate}
\end{prop}
\begin{pf}
(i) If $x \in L^{\perp} \cap (\cup_{\alpha \in R} H_{\alpha})$, 
$\langle x,l\rangle =0$ for any $l \in L$ and 
$\langle x,\alpha \rangle=0$ for some $\alpha \in R$. 
Since $\alpha \in \widetilde{F}$, 
$\alpha=l+Aa+B\delta$ for some $l \in L$ and $A,B \in \RR$. 
Then $0=A\langle a,x \rangle +B\langle \delta,x \rangle$. 
Since $\langle a,x \rangle,\ \langle \delta,x \rangle\in \CC$ are 
linearly independent over $\RR$ by $x \in L^{\perp}$, $A=B=0$. 
Then $\alpha \in L$. It contradicts the assumption of admissibility
 (\ref{200322.8}). 

(ii) Put $x^0:=z^0$. 
Since the set of zeros of determinant of 
$$
\left(
\frac{\p x^{\alpha}}{\p z^{\beta}}
\right)
_{0 \leq \alpha,\beta \leq n}
$$ 
coincides with 
$\displaystyle{\cup_{\alpha \in R}H_{\alpha}}$ 
(cf. \cite[(4.5) Theorem]{extendedII}), 
it is not $0$ on any point of $L^{\perp}$. 
By 
\begin{equation}
\frac{\p x^0}{\p z^0}=1,\quad
\frac{\p x^0}{\p z^{\alpha}}=0\quad (\one), 
\end{equation}
this determinant equals the determinant of 
$$
\left(
\frac{\p x^{\alpha}}{\p z^{\beta}}
\right)_{1 \leq \alpha,\beta \leq n}. 
$$
Then we have the result. 

(iii) By the result of (ii),  (\ref{200322.7}) and 
Proposition \ref{230402.1}, 
we obtain (iii) (where we used a discussion which is 
the same as the proof of Theorem 4.2(v) of \cite{Springer}). 
\qed\end{pf}
From now on we may and shall assume that 
\begin{equation}\label{200316.6}
g \cdot z^{\alpha}=\zeta^{d_{\alpha}}z^{\alpha}\quad (1 \leq \alpha \leq n).
\end{equation}
A $\CC$-basis of $L \otimes_{\RR}\CC$ satisfying (\ref{200316.6}) is called 
a ``$g$-homogeneous basis". 
\begin{prop}
We have 
\begin{eqnarray}
&&x^{\alpha}|_{L^{\perp}}=0\ (1 \leq \alpha \leq n,\ d_{\alpha}<d_n),
\label{200310.2}\\
&&\left.\frac{\p x^{\alpha}}{\p z^{\beta}}\right|_{L^{\perp}}=0\ 
(d_{\alpha}\neq d_{\beta}). 
\label{200316.1}
\end{eqnarray}
For any $a,b \in \multi$, we have 
\begin{equation}
\left.\frac{\p x^a}{\p z^b}\right|_{L^{\perp}}=0\ 
(a,b \in \multi,\ d\cdot b \not\equiv d\cdot a\ (\mathrm{mod}\ d_n)),
\label{200310.3}
\end{equation}
where we denote 
\begin{equation}
\frac{\p^b}{\p z^b}=
\left(\frac{\p}{\p z^1}\right)^{b_1}
\cdots
\left(\frac{\p}{\p z^n}\right)^{b_n}\quad
\mbox{ for }b=(b_1,\cdots,b_n) \in \multi. 
\end{equation}
\end{prop}
\begin{pf}
These are direct consequences of (\ref{200322.7}) and (\ref{200316.6}). 
\qed\end{pf}
%
%
%

%%%%%%%%%%%%%%%%%%%%%%%%%%%%%%%%%%%%%%%%%%%%%%%%%%%%%%%%%%%
%
%
%
%

%
%
%\newpage
%%%%%%%%%%%%%%%%%%%%%%%%%%%%%%%%%%%%%%%%%%%%%%%%%%%%%%%%%%

\section{Graded algebra Isomorphism $\psi$}
We fix an admissible triplet $(g,\zeta,L)$ for an elliptic root system. 
\subsection{The morphism $\varphi[g,\zeta,L]$}
\begin{defn}\label{210721.4}
For the admissible triplet $(g,\zeta,L)$, 
a set of basic invariants $x^1,\cdots,x^n$ and 
a $g$-homogeneous basis $z^1,\cdots,z^n$ of $L \otimes_{\RR}\CC$, 
we define an $F(H)$-module homomorphism:
\begin{equation}
\varphi[g,\zeta,L]: S^W \to F(Y),\quad
x^a \mapsto \sum_{b \in \multi}
\frac{1}{b!}
\left.\frac{\partial^b([x-(x|_{L^{\perp}})]^a)}{\partial z^b}\right|_{L^{\perp}}z^b,\label{200310.13}
\end{equation}
for an $F(H)$-free basis $\{x^a\,|\,a \in \multi\}$ 
of $F(H)$-module $S^W$, where we used notations
\begin{eqnarray}
[x-(x|_{L^{\perp}})]^a&:=&[x^1-(x^1|_{L^{\perp}})]^{a_1}\cdots 
[x^n-(x^n|_{L^{\perp}})]^{a_n} 
\ \mbox{for}\ (a=(a_1,\cdots,a_n) \in \ZZ^n_{\geq 0}),\nonumber \\
b!&:=&b_1!\cdots b_n!
\ \mbox{for}\ (b=(b_1,\cdots,b_n) \in \ZZ^n_{\geq 0}).
\end{eqnarray}
\end{defn}
We remark that $\varphi[g,\zeta,L](f)$ for $f \in S^W$ 
is not necessarily invariant by the $W$-action. 

\begin{prop}\label{200310.12}
\begin{enumerate}
\item $\varphi[g,\zeta,L]$ depends neither on the choices of 
a set of basic invariants $x^{\alpha}$ 
$(\one)$ nor on the choice of a set of $g$-homogeneous basis $z^{\alpha}$ $(\one)$ 
of $L \otimes_{\RR}\CC$. 
$\varphi[g,\zeta,L]$ gives an $F(H)$-algebra homomorphism. 
\item Let $x^{\alpha}\,(1 \leq \alpha \leq n)$ 
and $z^{\alpha}\,(1 \leq \alpha \leq n)$ be the same as 
in Definition $\ref{210721.4}$. 
For any multi-indices $a,b \in \multi$, the coefficients of $z^b$ 
of the RHS of 
$$
\varphi[g,\zeta,L](x^a)
=
\sum_{b \in \multi}
\frac{1}{b!}
\left.\frac{\partial^b([x-(x|_{L^{\perp}})]^a)}{\partial z^b}\right|_{L^{\perp}}z^b
$$
is $0$ if $d\cdot b \notin \{d\cdot a+d_n j\,|\,j \in \ZZ_{\geq 0}\}$. 
\end{enumerate}
\end{prop}
\begin{pf}
(i) For a set of basic invariants $x^{\alpha}\ (1 \leq \alpha \leq n)$, 
we define an $F(H)$-algebra homomorphism $\varphi_1[L]$ by 
\begin{equation}
\varphi_1[L]:S^W \to S^W,\quad 
x^{a} \mapsto [x-(x|_{L^{\perp}})]^a\ (a \in \multi). 
\end{equation}
By the same argument as in the proof of Proposition 3.2 (i) in 
\cite{handai3}, we see that this morphism does not depend 
on the choice of a set of basic invariants 
$x^1,\cdots,x^n$ but depends only on the choice of $L$. 

We define an $F(H)$-algebra homomorphism $\varphi_2$ by 
\begin{equation}
\varphi_2:S^W \to F(Y),\quad
f \mapsto \sum_{b \in \multi}\frac{1}{b!}
\left.\frac{\partial^b f}{\partial z^b}\right|_{L^{\perp}}
[z-(z|_{L^{\perp}})]^b. 
\end{equation}
This is a Taylor expansion along $L^{\perp}$ and it coincides 
with the natural inclusion 
$S^W \subset F(Y)$. 

We define an $F(H)$-algebra homomorphism $\varphi_3[L]$ by 
\begin{equation}
\varphi_3[L]:F(Y) \to F(Y), 
\quad
z^{a} \mapsto [z+(z|_{L^{\perp}})]^a\quad
(a \in \multi). 
\end{equation}
This morphism does not depend on the choice of a $g$-homogeneous basis 
$z^1,\cdots,z^n$ of $L \otimes_{\RR}\CC$ 
because a basis is unique up to linear transformations. 

Then we have 
\begin{equation}
\varphi[g,\zeta,L]=\varphi_3[L] \circ \varphi_2 \circ \varphi_1[L]. 
\end{equation}
Since $\varphi_1[L]$, $\varphi_2$ and $\varphi_3[L]$ are $F(H)$-algebra homomorphisms, 
their composite morphism $\varphi=\varphi[g,\zeta,L]$ is also 
an $F(H)$-algebra homomorphism. 

(ii) We remark that 
$\displaystyle{
\frac{\partial^b (x^{\alpha}|_{L^{\perp}})}{\partial z^b}=0
}$ if $b \neq 0$ 
by 
$\displaystyle{
x^{\alpha}|_{L^{\perp}} \in F(H)
}$. 
Then by the same argument as in the proof 
of Proposition 3.2\hskip0.2mm(ii) in \cite{handai3}, 
we have (ii). 
\qed
\end{pf}
%
%
%%%%%%%%%%%%%%%%%%%%%%%%%%%%%%%%%%%%%%%%%%%%%%%%%%%%
\subsection{The morphism $\psi[g,\zeta,L]$}

In this subsection, 
we construct a graded $F(H)$-algebra isomorphism 
$\psi[g,\zeta,L]$ by the same argument as in \S 3.2 in \cite{handai3}. 

We define decreasing filtrations on $S^W$ 
and $F(Y)$ by 
\begin{eqnarray*}
&&F^{m}(S^W):=
\bigoplus_{j \geq m}S^W_j
,\\
&&F^{m}(F(Y)):=\{\sum_{b \in \multi} c_b z^b\in F(Y)\,|\,
c_b \in F(H) (b \in \multi), 
c_b=0 \mbox { if } d\cdot b \leq m-1\}
\end{eqnarray*}
respectively for $m \in \ZZ_{\geq 0}$. Then $S^W$ and $F(Y)$ are filtered $F(H)$-algebras 
and $\varphi[g,\zeta,L]$ is a filtered $F(H)$-algebra homomorphism 
by Proposition \ref{200310.12}\hskip0.2mm(ii). 
\begin{defn}\label{230804.1}
\begin{enumerate}
\item 
Let $\mathrm{gr}_F\varphi[g,\zeta,L]$ be the 
graded $F(H)$-algebra homomorphism induced by 
a filtered $F(H)$-algebra homomorphism $\varphi[g,\zeta,L]$: 
\begin{equation}
\mathrm{gr}_F\varphi[g,\zeta,L]:\mathrm{gr}_F(S^W) \to \mathrm{gr}_F(F(Y)), 
\end{equation}
where 
\begin{eqnarray}
\mathrm{gr}_F(S^W)
&:=&\bigoplus_{m \in \ZZ_{\geq 0}}F^{m}(S^W)/F^{m+1}(S^W),\\
\mathrm{gr}_F(F(Y))
&:=&\bigoplus_{m \in \ZZ_{\geq 0}}F^{m}(F(Y))/F^{m+1}(F(Y)). 
\end{eqnarray}
\item For the graded $F(H)$-algebra $S^W$, 
we have the natural graded $F(H)$-algebra isomorphism 
\begin{equation}
\psi_1:S^W \to \mathrm{gr}_F(S^W)
\end{equation}
which maps an element of $S^W_j$ to its canonical image in 
$F^j(S^W)/F^{j+1}(S^W)$. 
Let $z^1,\cdots,z^n$ be a $g$-homogeneous basis of $L \otimes_{\RR}\CC$. 
We define 
\begin{equation}
F(H)[L]:=\{\sum_{b \in \multi} 
c_b z^b\in F(Y)\,|\,c_b \in F(H),\ d\cdot b \mbox{ is bounded }\}.
\end{equation}
We have a decomposition
\begin{equation}\label{200309.2}
F(H)[L]=\bigoplus_{j \in \ZZ}V(g,\zeta,L)(j), 
\end{equation}
where 
\begin{equation}\label{200316.7}
V(g,\zeta,L)(j):=
\{\sum_{b \in \multi} c_b z^b\in F(Y)\,|\,c_b \in F(H), d\cdot b=j\}.
\end{equation}
for $j \in \mathbb{Z}_{\geq 0}$. 
These definitions do not depend on the choice of 
$z^1,\cdots,z^n$. 
The decomposition (\ref{200309.2}) 
gives a graded $F(H)$-algebra structure on $F(H)[L]$ 
which is isomorphic to the polynomial algebra
\begin{equation}
F(H)[z^1,\cdots,z^n]
\end{equation}
with $\mathrm{deg}z^{\alpha}=d_{\alpha}$. 
Since the composite mapping 
\begin{equation}
V(g,\zeta,L)(j) \to F^j(F(Y)) \to F^j(F(Y))/F^{j+1}(F(Y))
\end{equation}
is an isomorphism, we have a graded $F(H)$-algebra isomorphism
\begin{equation}
\psi_2:F(H)(L) \to \mathrm{gr}_F(F(Y))). 
\end{equation}
\item Let $\psi[g,\zeta,L]$ be the graded $F(H)$-algebra homomorphism 
defined by 
\begin{equation}
\psi[g,\zeta,q]:=\psi_2^{-1} \circ \mathrm{gr}_F\varphi[g,\zeta,L] 
\circ \psi_1:
S^W \to F(H)(L). 
\end{equation}
\end{enumerate}
\end{defn}
We have an explicit description of $\psi[g,\zeta,L]$:
\begin{equation}
\psi[g,\zeta,L]:S^W \to F(H)[L],\quad
x^a \mapsto \sum_{b \in \multi,\ 
d\cdot b=d\cdot a}
\frac{1}{b!}
\left.\frac{\partial^b([x-(x|_{L^{\perp}})]^a)}{\partial z^b}\right|_{L^{\perp}}z^b\label{200310.14}
\end{equation}
for an $F(H)$-free basis $\{x^a\,|\,a \in \multi\}$ 
of $F(H)$-module $S^W$, where we used notations in Definition \ref{210721.4}.

By Proposition \ref{200322.12}\hskip0.2mm(ii), 
we have the following proposition 
by the same argument as in the proof of Proposition 3.3 in \cite{handai3}. 
\begin{prop}
With respect to the gradings $(\ref{200309.1})$ on $S^W$ 
and $(\ref{200309.2})$ on $F(H)[L]$, 
$\psi[g,\zeta,L]$ is a graded $F(H)$-algebra isomorphism
\begin{equation}
\psi[g,\zeta,L]:S^W
\stackrel{\sim}{\to}
F(H)[L].
\end{equation}
\end{prop}

\subsection{Good basic invariants}
\begin{defn}\label{200316.9}
A set of basic invariants $x^1,\cdots,x^n$ is good with respect to 
the admissible triplet $(g,\zeta,L)$ if 
$x^1,\cdots,x^n$ form a $\mathbb{C}$-basis of the vector space 
$\psi[g,\zeta,L]^{-1}(L\otimes_{\RR}\CC)$ 
w.r.t. the natural inclusion $L \otimes_{\RR}\CC \subset F(H)[L]$. 
We call $x^1,\cdots,x^n$ ``good basic invariants". \\
\end{defn}
%\newpage
\vskip1cm
%\newpage

\section{Taylor coefficients of the good basic invariants}
Let $(g,\zeta,L)$ be an admissible triplet 
for an elliptic root system. 
\begin{defn}\label{200316.10}
Let $z^0=\frac{\delta}{-2\pi\sqrt{-1}}$ defined in (\ref{200323.12}) 
and $z^1,\cdots,z^n$ be a $g$-homogeneous basis of $L \otimes_{\RR}\CC$. 
Then a set of basic invariants $x^1,\cdots,x^n$ is 
compatible with a basis $z^1,\cdots,z^n$ of $L\otimes_{\RR}\CC$ if 
the Jacobian matrix is a unit matrix, i.e. 
$$
\left(
\left.\frac{\p x^{\alpha}}{\p z^{\beta}}\right|_{L^{\perp}}
\right)_{1 \leq \alpha,{\beta} \leq n}
=\left(
\delta^{\alpha}_{\beta}
\right)_{1 \leq \alpha,{\beta} \leq n},
$$
where $\delta^{\alpha}_{\beta}$ is the Kronecker's delta. 
\end{defn}
\begin{prop}\label{230817.3}
For a $g$-homogeneous basis $z^1,\cdots,z^n$ of $L \otimes_{\RR}\CC$, 
we have the following results.
\begin{enumerate}
\item If we put 
\begin{equation}
x^{\alpha}:=\psi[g,\zeta,L]^{-1}(z^{\alpha})\quad
(\one),
\end{equation}
then $x^1,\cdots,x^n$ form a set of basic invariants 
which are good and compatible with a $g$-homogeneous 
basis $z^1,\cdots,z^n$ of $L\otimes_{\RR}\CC$. 
\item Conversely if $x^1,\cdots,x^n$ are good and 
compatible with a $g$-homogeneous basis $z^1,\cdots,z^n$ of $L\otimes_{\RR}\CC$, 
then $\psi[g,\zeta,L](x^{\alpha})=z^{\alpha}$ for $\one$. 
\item For any set of basic invariants $x^1,\cdots,x^n$, 
\begin{equation}\label{200310.9}
\left.\frac{\partial^b([x-x|\rest]^a)}{\partial z^b}\right|_{L^{\perp}}=0
\hbox{ if }
d \cdot b \notin \{d\cdot a+d_n j\,|\, j \in \mathbb{Z}_{\geq 0}\}
\end{equation}
for $a,b \in \multi$. 
\item A set of basic invariants $x^1,\cdots,x^n$ is good if and only if 
$$
\left.\frac{\p x^{\alpha}}{\p z^{\beta}}\right|_{L^{\perp}}
\quad
(1 \leq \alpha,{\beta} \leq n)
$$
are constants and 
\begin{equation}\label{200318.7}
\left.\frac{\p^a x^{\alpha}}{\p z^a}\right|_{L^{\perp}}=0\ 
(d_{\alpha}=d\cdot a,\ |a|\geq 2,\ \one). 
\end{equation}
%
%
%\item $x^{\alpha}$ are very good iff
%\begin{equation}
%\frac{\p^a x^{\alpha}}{\p z^a}(q)=0\ 
%(d_i=d\cdot a,\ |a|\geq 2),\quad
%\frac{\p^a x^{\alpha}}{\p z^{\beta}}(q)=\delta_{\alpha,\beta}.  
%\end{equation}
%
%
\item If a set of basic invariants $x^1,\cdots,x^n$ is good and compatible 
with a $g$-homogeneous basis $z^1,\cdots,z^n$ of $L\otimes_{\RR}\CC$, then 
for $a,b \in \multi$ satisfying $d \cdot a=d\cdot b$, we have 
\begin{equation}\label{200310.10}
\left.\frac{1}{b!}\frac{\partial^b([x-x|\rest]^a)}{\partial z^b}\right|_{L^{\perp}}=\delta_{a,b}.
\end{equation}
\item If a set of basic invariants $x^1,\cdots,x^n$ is good, 
then for $1 \leq \alpha \leq n$ and $a \in \multi$ 
satisfying $d \cdot a=d_{\alpha}$, we have 
\begin{equation}\label{200323.13}
\left.\left(
\frac{\p}{\p z^0}
\frac{1}{a!}
\frac{\partial^{a}x^{\alpha}}
{\partial z^a}
\right)
\right|_{L^{\perp}}=0.
\end{equation}
\end{enumerate}
\end{prop}
\begin{pf} As for (i), (ii), they are direct consequences of 
Definition \ref{200316.9}
and 
Definition \ref{200316.10}. 

As for (iii), it is proved in Proposition \ref{200310.12} (ii). 
As for (iv), we have 
$$
\psi[g,\zeta,L](x^{\alpha})=
\sum_{b \in \multi,\ 
d\cdot b=d_{\alpha}}\frac{1}{b!}
\left.\frac{\partial^b x^{\alpha}}{\partial z^b}\right|_{L^{\perp}} z^b
$$
by (\ref{200310.14}). 
By the goodness assumption, this must be an element of $L\otimes_{\RR}\CC$. 
Then the coefficient of $z^b$ is constant if $|b|=1$ 
and $0$ if $|b| \geq 2$. 

As for (v), we have $\psi[g,\zeta,L](x^{\alpha})=z^{\alpha}$ for $\one$ 
by (ii). Then for any $a \in \multi$, 
\begin{equation}
\psi[g,\zeta,L](x^a)
=\prod_{\gamma=1}^n \psi[g,\zeta,L](x^{\gamma})^{a_i}
=\prod_{\gamma=1}^n (z^{\gamma})^{a_i}
=z^a, 
\end{equation}
and comparing it with (\ref{200310.14}), we have the result. 
As for (vi), we have 
$$
\left.
\left(
\frac{\p}{\p z^0}
\frac{1}{a!}
\frac{\partial^{a}x^{\alpha}}
{\partial z^a}
\right)
\right|_{L^{\perp}}
=
\frac{\p}{\p z^0}
\left(
\left.
\frac{1}{a!}
\frac{\partial^{a}x^{\alpha}}
{\partial z^a}\right|_{L^{\perp}}
\right)
$$
and 
$$
\left.
\frac{1}{a!}
\frac{\partial^{a}x^{\alpha}}
{\partial z^a}\right|_{L^{\perp}}
$$
is constant because $x^{\alpha}$ is good. Then we have the result. 
\qed
\end{pf}

\section{Construction of an admissible triplet}
\subsection{Construction of an admissible triplet}
In this subsection, we construct an admissible triplet 
for an elliptic root system. 

The following theorem is due to Saito \cite{extendedI}. 
\begin{thm}\label{200324.1}
$($Saito \cite{extendedI}$)$ 
There exists $\widetilde{c} \in W \subset O(\widetilde{F},F,\mathrm{rad}\,I)$ 
called a hyperbolic Coxeter transformation $($\cite[(11.2)]{extendedI}$)$ 
which satisfies the following properties.
\begin{enumerate}
\item $($Lemma A$)$ The restriction $c$ of $\widetilde{c}$ to $F$ 
$($called a Coxeter transformation \cite[(9.7)]{extendedI}$)$ 
is semi-simple of order $d_n$. 
The set of eigenvalues of $c$ is given by:
\begin{equation}
1,\ 
\exp
\left(
2\pi\sqrt{-1}d_{\alpha}/d_n
\right)
\quad
(\alpha=1,\cdots,n).
\end{equation}
\item $($Lemma B$)$ Let $c$ be a Coxeter transformation. Then 
\begin{equation}
R \cap \mathrm{Im}(c-id.)=\emptyset.
\end{equation}
\item $($Lemma C$)$ 
For a hyperbolic Coxeter transformation $\widetilde{c}$, we have 
\begin{equation}
(\widetilde{c}-1)\xi+I_R(\xi,\delta)\frac{1}{m_{max}}a \in \mathrm{Im}(c-id.)
\quad
(\forall \xi \in \widetilde{F})
\end{equation}
and $\widetilde{c}^{d_n}$ is a generator of $K_{\ZZ}$, 
where $I_R=(I_R:I)I$, $(I_R:I)$ is defined in Section \ref{240311.2} 
and $m_{max}$ is defined in Section \ref{240.001}.
\end{enumerate}
\end{thm}
Let 
\begin{equation}
\widetilde{c}=\widetilde{c}^{ss}\cdot \widetilde{c}^{unip}
\end{equation}
be the Jordan decomposition to semi-simple element and 
unipotent element. 

The following proposition shows 
that the semi-simple element $\widetilde{c}^{ss}$ satisfies 
(\ref{230327.2}) and (\ref{200322.7}) which are part of 
conditions of the admissible triplet. 
\begin{prop}\label{200322.5}
\begin{enumerate}
\item $\widetilde{c}^{ss}$ is an element of $O(\widetilde{F},F,\mathrm{rad}\,I)$ 
and 
\begin{equation}
\widetilde{c}^{ss}|_F=c.
\end{equation}
\item 
$$
(\widetilde{c}^{ss}-id.)(\widetilde{F}) \subset F.
$$
\item There exists uniquely a primitive $d_n$-th root of unity $\zeta$ 
such that the action of $\widetilde{c}^{ss}$ on $f \in S^W_m$ is given by 
\begin{equation}
\widetilde{c}^{ss}\cdot f=\zeta^{m} f.
\end{equation}
\end{enumerate}
\end{prop}
\begin{pf}
We first show that the $\RR$-bilinear form $I$ on 
$(c-id.)(F)$ is positive definite. 
Since $c-id.:F \to F$ is semi-simple, 
we have a decomposition:
$$
F=(c-id.)(F) \oplus \mathrm{ker}(c-id.).
$$
By the inclusion $\mathrm{rad}\,I \subset \mathrm{ker}(c-id.)$, 
we have $(c-id.)(F) \cap \mathrm{rad}\,I=\emptyset$. 
Since $I$ on $F$ is semi-positive, $I$ on $(c-id.)(F)$ is 
positive definite. 

Put 
$$
(c-id.)(F)^{\perp}:=\{x \in \widetilde{F}\,|\,
\widetilde{I}(x,y)=0\ 
\forall y \in (c-id.)(F)\}.
$$
Then we have 
\begin{equation}\label{230327.3}
\widetilde{F}=(c-id.)(F) \oplus (c-id.)(F)^{\perp}.
\end{equation}
We see that $(c-id.)(F)$ is $\widetilde{c}$-stable and 
$\widetilde{c}$ on $(c-id.)(F)$ is semi-simple because 
$\widetilde{c}=c$ on $F$ and $c$ is semi-simple. 

The space $(c-id.)(F)^{\perp}$ is $\widetilde{c}$-stable 
because $(c-id.)(F)$ is $\widetilde{c}$-stable and 
$\widetilde{c} \in O(\widetilde{F},F,\mathrm{rad}\,I)$. 

We show that the action of $\widetilde{c}$ on 
$(c-id.)(F)^{\perp}$ is unipotent. 

For all $x \in (c-id.)(F)^{\perp}$, 
\begin{equation}\label{230327.4}
(\widetilde{c}-id.)(x)+I_R(x,\delta)\frac{1}{m_{max}}a
\end{equation}
is an element of $(c-id.)(F)$ by Theorem \ref{200324.1}
(iii). 
On the other hand, $(\widetilde{c}-id.)(x)$ in (\ref{230327.4}) 
is an element of $(c-id.)(F)^{\perp}$ because 
$(c-id.)(F)^{\perp}$ is $\widetilde{c}$-stable. 
Also $\displaystyle{I_R(x,\delta)\frac{1}{m_{max}}a}$ in 
(\ref{230327.4}) is an element of $\mathrm{rad}\,I$ 
which is a subset of $(c-id.)(F)^{\perp}$. 
Then by a decomposition (\ref{230327.3}), we have 
\begin{equation}\label{230327.5}
(\widetilde{c}-id.)(x)+I_R(x,\delta)\frac{1}{m_{max}}a=0.
\end{equation}
Then $\widetilde{c}$ is identity both on $(c-id.)(F)^{\perp}/\RR a$ 
and $\RR a$. 
Thus $\widetilde{c}$ is unipotent on $(c-id.)(F)^{\perp}$. 

Then we have 
\begin{eqnarray}
\widetilde{c}^{ss}=
\begin{cases}
\widetilde{c} &\mbox{ on }(c-id.)(F)\\
id. &\mbox{ on }(c-id.)(F)^{\perp},
\end{cases}
\quad
\widetilde{c}^{unip}=
\begin{cases}
id. &\mbox{ on }(c-id.)(F)\\
\widetilde{c} &\mbox{ on }(c-id.)(F)^{\perp}.
\end{cases}
\end{eqnarray}
From this description of $\widetilde{c}^{ss}$, 
we have (i), (ii). 

As for (iii), we see that 
$\widetilde{c}^{unip}$ is an element of $K_{\RR}$ (see (\ref{200322.4})) 
and $(\widetilde{c}^{unip})^{d_n}$ is a generator of $K_{\ZZ}$ 
because $(\widetilde{c}^{unip})^{d_n}$ is a unipotent part of 
$\widetilde{c}^{d_n}$ and it is a generator of $K_{\ZZ}$ by Theorem 
\ref{200324.1}\hskip0.2mm(iii). 

Then there exists uniquely a primitive $d_n$-th root of unity $\zeta$ 
such that the action of $\widetilde{c}^{unip}$ on $f \in S^W_m$ is given by 
\begin{equation}
\widetilde{c}^{unip}\cdot f=\zeta^{-m} f. 
\end{equation}
Since $\widetilde{c}^{ss}=(\widetilde{c}^{unip})^{-1}\cdot \widetilde{c}$ 
and $\widetilde{c} \in W$, the action of $\widetilde{c}^{ss}$ on $f \in S^W_m$ is 
$$
\widetilde{c}^{ss}\cdot f
=
[(\widetilde{c}^{unip})^{-1}\cdot \widetilde{c}]\cdot f
=
(\widetilde{c}^{unip})^{-1}\cdot f
=
\zeta^{m}f.
$$
\qed\end{pf}
%\newpage
We shall construct a splitting subspace $L$ 
which is a part of admissible triplet. 
We remind the reader of the definition of the space $\widetilde{X}$ defined 
in $(\ref{230812.5})$:
$$
\widetilde{X}=\{V \subset \widetilde{F}\,|\,
\widetilde{F}=V \oplus \mathrm{rad}\,I\}.
$$
\begin{defn}
We put 
%a set of regular subspaces of $\widetilde{F}$ 
%and a set of fixed subspaces of $\widetilde{F}$ w.r.t. 
%the action of $\widetilde{c}^{ss}$:
\begin{eqnarray}
\widetilde{X}^{reg}&:=&
\{V \in \widetilde{X}\,|\,V \cap R= \emptyset\},\\
\widetilde{X}^{\widetilde{c}^{ss}}&:=&
\{V \in \widetilde{X}\,|\,\widetilde{c}^{ss}\cdot V =V\}.
\end{eqnarray}
\end{defn}

\begin{prop} \label{240227.1}
We have 
\begin{equation}
\widetilde{X}^{\widetilde{c}^{ss}} \cap \widetilde{X}^{reg}\neq \emptyset.
\end{equation}
\end{prop}
A proof of Proposition \ref{240227.1} will be given in the next subsection. 
\begin{prop}\label{240228.1}
Let $\zeta$ be a primitive $d_n$-th root of unity defined 
in Proposition $\ref{200322.5}\ (iii)$. 
For any $L \in \widetilde{X}^{\widetilde{c}^{ss}} \cap \widetilde{X}^{reg}$, 
$(\widetilde{c}^{ss},\zeta,L)$ is an admissible triplet. 
\end{prop}
\begin{pf}
By Proposition \ref{200322.5}, $\widetilde{c}^{ss}$ and $\zeta$ 
satisfies the conditions in Definition \ref{200310.1}. 
For $L \in \widetilde{X}^{\widetilde{c}^{ss}} \cap \widetilde{X}^{reg}$, 
it satisfies the condition of Definition \ref{200310.1}\hskip0.2mm(iii). 
\qed
\end{pf}

\subsection{Proof of Proposition \ref{240227.1}}
In this subsection, we give a proof of Proposition \ref{240227.1}. 

We first prepare a space of complementary subspaces of 
$\mathrm{rad}\,I$ in a vector space $F$:
\begin{equation}
X:=\{U \subset F\,|\,F=U \oplus \mathrm{rad}\,I\}. 
\end{equation}
We give a relation of $X$ with $\widetilde{X}$ 
defined in $(\ref{230812.5})$. 
For $V \in \widetilde{X}$, $V \not\subset F$, 
the dimension of $V \cap F$ is $n-1(=l)$. 
Then we have 
$$
F=(V \cap F) \oplus \mathrm{rad}\,I. 
$$
Then we have a natural morphism:
\begin{equation}
p:\widetilde{X} \to X,\quad
V \mapsto V \cap F.
\end{equation}
We could easily check that 
$p:\widetilde{X} \to X$ is a 
$\mathrm{Hom}_{\RR}(\widetilde{F}/F,\mathrm{rad}\,I)$-principal 
bundle. 
In particular, $p$ is surjective. 
The group $O(\widetilde{F},F,\mathrm{rad}\,I)$ acts on the 
space $X$ and $p$ is equivariant w.r.t. this group action. 
\begin{defn}
We put 
%a set of regular subspaces of $F$ 
%and a set of fixed subspaces of $F$ w.r.t. 
%the action of $\widetilde{c}^{ss}$:
\begin{eqnarray}
X^{reg}&:=&
\{U \in X\,|\,U \cap R= \emptyset\},\\
X^{\widetilde{c}^{ss}}&:=&
\{U \in X\,|\,\widetilde{c}^{ss}\cdot U =U\}.
\end{eqnarray}
\end{defn}
\begin{prop}\label{240227.3}
\begin{eqnarray}
\mathrm{(i)}\ 
p^{-1}(X^{reg})=\widetilde{X}^{reg}.
\label{230330.10}
\\
\mathrm{(ii)}\ 
p^{-1}(X^{\widetilde{c}^{ss}})=\widetilde{X}^{\widetilde{c}^{ss}}.
\label{230330.9}
\end{eqnarray}
\end{prop}
\begin{pf}
(i) is a consequence of Lemma \ref{230329.3}. 
(ii) is a consequence of 
Proposition \ref{200322.5}\hskip0.2mm(ii) 
and 
Lemma \ref{230330.6}\hskip0.2mm(iii). 
\qed
\end{pf}
\begin{lem}\label{230329.3}
For $V \in \widetilde{X}$, 
$V \cap R =\emptyset$ if and only if 
$p(V) \cap R=\emptyset$. 
\end{lem}
\begin{pf} If $V \cap R =\emptyset$, $(V \cap R) \cap F$ is also 
an empty set, which is $p(V) \cap R$. 

If $V \cap R \neq \emptyset$, 
there exists $x \in V \cap R$. Then $x \in R \subset F$, 
$(V \cap R) \cap F$ is also non-empty. 
\qed
\end{pf}
\begin{lem} \label{230330.6}
Let $g \in O(\widetilde{F},F,\mathrm{rad}\,I)$ which is semi-simple
with $(g-id.)(\widetilde{F}) \subset F$. 
\begin{enumerate}
\item We have
$$
(g-id.)(F)=(g-id.)(\widetilde{F}).
$$
\item For $V \in \widetilde{X}$, the followings are equivalent.
\begin{enumerate}
\item $g(V)=V$.
\item $g(p(V))=p(V)$. 
\item $(g-id.)(\widetilde{F}) \subset V$. 
\end{enumerate}
\item 
We put
% a set of fixed points:
$$
X^g:=\{U \in X\,|\,g\cdot U=U\}. 
$$
Then we have 
$$
p^{-1}(X^g)=
\{V \in \widetilde{X}\,|\,
g(V)=V\}.
$$
\end{enumerate}
\end{lem}
\begin{pf} 
(i) The proof of the inclusion $\subset$ is trivial. For $\supset$, 
we have $(g-id.)^2(\widetilde{F}) \subset (g-id.)(F)$ by the assumption. 
Since $g-id.:\widetilde{F} \to \widetilde{F}$ is semi-simple, 
$(g-id.)^2(\widetilde{F})=(g-id.)(\widetilde{F})$ Then we have a result. 

(ii) 
(a) $\Longrightarrow$ (b). 
Since $g(F) \subset F$, 
$g(F \cap V) \subset F \cap V$. 

(b) $\Longrightarrow$ (c). 
Since $p(V)$ is $g$-stable, 
a natural projection 
$\varphi:F \to F/p(V) \simeq \mathrm{rad}\,I$ 
is $g$-equivalent. 
For any $x \in (g-id.)(F)$, $\exists y \in F$ s.t. 
$x=(g-id.)(y)$. Then we have 
$$
\varphi(x)=\varphi(g(x)-x)=g\cdot \varphi(x)-\varphi(x)=0, 
$$
because the action of $g$ on $\mathrm{rad}\,I$ is trivial. 
Then $(g-id.)(F) \subset p(V)$. 

By $(g-id.)(\widetilde{F})=(g-id.)(F)$ and 
$p(V) \subset V$, we have (c). 

(c) $\Longrightarrow$ (a). 
Since $g-id.:\widetilde{F} \to \widetilde{F}$ is semi-simple, 
a vector space $V$ which satisfies 
$$
(g-id.)(F) \subset V \subset \widetilde{F}
$$
is $g$-stable. 

(iii) is a direct consequence of (ii).
\qed
\end{pf}

For the space and a natural morphism:
\begin{eqnarray}
&\mathbb{E}:=
\{
x \in \mathrm{Hom}_{\RR}(F,\CC)\,|\,\langle a,x \rangle=-2\pi\sqrt{-1},\ 
\mathrm{Re} \langle \delta,x \rangle>0\,
\},\label{230330.4}\\
&\pi':\mathbb{E} \to H,\label{230330.5}
\end{eqnarray}
we have a similar construction as in \S \ref{230813.4}, i.e. 
we define a mapping:
\begin{equation}
f_2:\mathbb{E} \to X
\end{equation}
by $f_2(x)=\mathrm{ker}\ x$ where we see $x \in \mathbb{E}$ 
as a morphism $x:F \to \CC$. 
We see that $f_2$ is $O(F,\mathrm{rad}\,I)$-equivariant. 
We have a mapping:
\begin{equation}\label{230330.7}
(\pi',f_2):\mathbb{E} \to H \times X
\end{equation}
which is an isomorphism as a real manifold. 

The following proposition which is obtained 
by Theorem \ref{200324.1} (ii) Lemma B 
is due to Saito \cite{extendedII} 
\begin{prop}$($\cite[p.44 (7.2) Lemma]{extendedII}$)$ \label{240227.2}
For a Coxeter transformation in Theorem $\ref{200324.1}(\mathrm{i})$, 
we put 
\begin{eqnarray}
&&\mathbb{E}^{c}:=\{x \in \mathbb{E}\,|\,
c\cdot x=x\ \},\\
&&\mathbb{E}_{\tau}:=(\pi')^{-1}(\tau)
\end{eqnarray}
for the space $\mathbb{E}$ defined in $(\ref{230330.4})$, 
$\pi':\mathbb{E} \to H$ defined in $(\ref{230330.5})$ and $\tau \in H$. 
Then we have 
\begin{equation}
\mathbb{E}_{\tau} \cap \mathbb{E}^c \not\subset 
\bigcup_{\alpha \in R}\overline{H}_{\alpha}
\end{equation}
for 
\begin{eqnarray}
\overline{H}_{\alpha}:=\{x \in \mathbb{E}\,|\,
\langle \alpha,x\rangle=0\ \}\quad (\alpha \in R).
\end{eqnarray}
\end{prop}
\begin{prop} \label{230330.11}
We have 
\begin{equation}\label{230402.3}
X^{\widetilde{c}^{ss}} \cap X^{reg}\neq \emptyset.
\end{equation}
\end{prop}
\begin{pf}
For $\tau \in H$, 
the isomorphism (\ref{230330.7}) induces an isomorphism:
\begin{equation}\label{230330.8}
\mathbb{E}_{\tau} \simeq X.
\end{equation}
We could easily check that this isomorphism 
induces the following isomorphisms:
\begin{eqnarray}
&&\mathbb{E}_{\tau} \cap \mathbb{E}^c \simeq X^{\widetilde{c}^{ss}},
\label{240309.4}\\
&&\mathbb{E}_{\tau}\setminus \bigcup_{\alpha \in R}\overline{H}_{\alpha}
\simeq X^{reg}
\end{eqnarray}
by Proposition \ref{200322.5}\hskip0.2mm(i). 
Then by Proposition \ref{200324.1}\hskip0.2mm(iii) and Proposition \ref{240227.2}, we have (\ref{230402.3}). 

\qed
\end{pf}

Proof of Proposition \ref{240227.1}. 
By Proposition \ref{230330.11}, we have 
$X^{\widetilde{c}^{ss}} \cap X^{reg}\neq \emptyset$. 
Since $p:\widetilde{X} \to X$ is surjective, 
we have $p^{-1}(X^{\widetilde{c}^{ss}} \cap X^{reg})\neq \emptyset$. 
By Proposition \ref{240227.3}, we see that 
$$
\widetilde{X}^{\widetilde{c}^{ss}}\cap \widetilde{X}^{reg}
=
p^{-1}(X^{\widetilde{c}^{ss}})\cap p^{-1}(X^{reg})
=
p^{-1}(X^{\widetilde{c}^{ss}} \cap X^{reg})
\neq \emptyset.
$$ 
\qed
\begin{rmk}
The subspace $L$ which we construct in Proposition \ref{240228.1} 
satisfies $(c-id.)(F) \subset L$ by Lemma \ref{230330.6}\hskip0.2mm(i), (ii). 
Then the assertion $R \cap L=\emptyset$ is 
an enhancement of Theorem \ref{200324.1} (ii) Lemma B. 
\end{rmk}

\section{Ambiguity of the choice of a splitting subspace $L$}\label{240314.2}
In this section, we show that the $\CC$-span of good basic invariants 
do not depend on the choice of a splitting subspace 
$L \in \widetilde{X}^{\widetilde{c}^{ss}} \cap \widetilde{X}^{reg}$ 
of an admissible triplet $(\widetilde{c}^{ss},\zeta,L)$ 
of ``zero type" which we define in Definition \ref{240309.1} 
if the codimension of an elliptic root system is $1$. 

For an elliptic root system, 
its codimension is defined in \cite[p23]{extendedII} as a cardinarity of 
$\{i|d_i=d_n\}$. 

If the codimension of an elliptic root system is $1$, 
then we have a uniqueness assertion (Proposition \ref{240309.2}) 
that the $\CC$-span of good basic invariants 
for the admissible triplet $(\widetilde{c}^{ss},\zeta,L)$ 
does not depend on the choice of 
$L \in \widetilde{X}^{\widetilde{c}^{ss}} \cap \widetilde{X}^{reg}$ 
under the assumption that $L$ is of ``zero type" 
which we define in Definition \ref{240309.1}. 

If the codimension of an elliptic root system is greater than $1$, 
then we have no uniqueness theorem. 
In Appendix A, we give an example of 
$L_1,L_2 \in \widetilde{X}^{\widetilde{c}^{ss}} \cap \widetilde{X}^{reg}$ 
which are of zero type and 
give different $\CC$-spans of good basic invariants 
for the case of $A_1^{(1,1)}$ type. 

\subsection{Signature of $V \in \widetilde{X}$}
For $V \in \widetilde{X}$, 
$\widetilde{I}$ on $V \cap F$ is positive definite. 
Then a signature of $\widetilde{I}$ on $V$ may be 
$(n-1,1,0)$ or $(n,0,0)$ or $(n-1,0,1)$, 
where we denote by $(l_{+},l_{0},l_{-})$ the numbers 
of positive, zero and negative eigenvalues 
of the Gram matrix of $\widetilde{I}|_V$.
\begin{defn} \label{240309.1}
An element $V \in \widetilde{X}$ is called of ``zero type" if 
a signature of $\widetilde{I}$ on $V$ is $(n-1,1,0)$,  
\end{defn}
We give an explicit description of a splitting subspace of zero type. 
$V \in \widetilde{X}$. 
We remind that $a,\delta \in \mathrm{rad}\,I$ is a basis of 
$\mathrm{rad}\,I$. 
For $U \in X$, an $\RR$-vector space
$$
\{ x \in \widetilde{F}\,|\,
\widetilde{I}(x,y)=0\ \forall y \in U\}
$$
has an $\RR$-basis $a,\delta,\lambda$ such that 
$$
\widetilde{I}(\delta,\lambda)=1,\quad
\widetilde{I}(\lambda,\lambda)=0.
$$
For $c_1,c_2 \in \RR$, put
\begin{equation}\label{240309.3}
V_{c_1,c_2}:=U \oplus \RR(\lambda+c_1\delta+c_2 a).
\end{equation}
For $p:\widetilde{X} \to X$, we have
$$
p^{-1}(U)=\{V_{c_1,c_2}\,|\,c_1,c_2 \in \RR\}.
$$
Then we see that $V_{c_1,c_2}$ is of zero type 
if and only if $c_1=0$ because $\widetilde{I}$ on $U$ is 
positive definite. 

\subsection{Uniqueness of the good basic invariants}
\begin{prop}\label{240309.5}
Let $x^1,\cdots,x^n$ be a set of good basic invariants 
for the admissible triplet 
$(\widetilde{c}^{ss},\zeta,V_{0,0})$, 
where $V_{0,0}$ is constructed in $($\ref{240309.3}$)$ 
for $U \in X^{\widetilde{c}^{ss}} \cap X^{reg}$. 
Then for $c_1,c_2 \in \RR$, 
$$
\widetilde{x}^1=(e^{c_1 \delta+c_2 a})^{d_1}x^1
,\cdots,
\widetilde{x}^n=(e^{c_1 \delta+c_2 a})^{d_n}x^n
$$
are a set of good basic invariants for the admissible 
triplet $(\widetilde{c}^{ss},\zeta,V_{c_1,c_2})$. 
\end{prop}
\begin{pf}
Let $z^0=\frac{\delta}{-2\pi\sqrt{-1}}$ defined in (\ref{200323.12}). 
We take a $g$-homogeneous basis $z^1,\cdots,z^{n-1} \in U \otimes_{\RR}\CC$. 
We take $z^n \in V_{0,0}$ such that $z^1,\cdots,z^n$ 
is a $g$-homogeneous basis of $V \otimes_{\RR}\CC$ and 
$E z^n=1$, where $E$ is a Euler field defined in (\ref{200322.1}). 
Put 
$$
\widetilde{z}^0:=z^0,\ 
\widetilde{z}^1:=z^1,\cdots,
\widetilde{z}^{n-1}:=z^{n-1},\ 
\widetilde{z}^n:=z^n+c_1\delta+c_2 a.
$$
Then a set of $\widetilde{z}^1,\cdots,\widetilde{z}^{n}$ is a $g$-homogeneous 
basis of $V_{c_1,c_2} \otimes_{\RR}\CC$. 

If $f \in F(Y)$ satisfies $Ef=mf$ for $m \in \ZZ$, 
then $f$ has a decomposition:
$$
f=\exp(mz^n)\overline{f}(z^0,z^1,\cdots,z^{n-1})
$$
for some function $\overline{f}(z^0,z^1,\cdots,z^{n-1})$. 
Then 
$$
f=\exp(m(\widetilde{z}^n-c_1\delta-c_2 a))\overline{f}(\widetilde{z}^0,
\widetilde{z}^1,\cdots,\widetilde{z}^{n-1}).
$$
Thus 
$$
f|_{V_{0,0}^{\perp}}=\overline{f}(z^0,0,\cdots,0)
$$
and
$$
f|_{V_{c_1,c_2}^{\perp}}=
\exp(m(-c_1\delta-c_2 a))\overline{f}(\widetilde{z}^0,
0,\cdots,0).
$$
Then we have 
\begin{equation}\label{230817.2}
f|_{V_{c_1,c_2}^{\perp}}
=\exp(m(-c_1\delta-c_2 a))
f|_{V_{0,0}^{\perp}},
\end{equation}
where we compare these functions under the identification 
$V_{0,0}^{\perp}\simeq H \simeq V_{c_1,c_2}^{\perp}$. 

We have 
\begin{eqnarray*}
&&\left.\frac{\p^a \widetilde{x}^{\alpha}}{\p \widetilde{z}^a}
\right|_{V_{c_1,c_2}^{\perp}}
=
(e^{c_1 \delta+c_2 a})^{d_{\alpha}}
\left.\frac{\p^a x^{\alpha}}{\p \widetilde{z}^a}
\right|_{V_{c_1,c_2}^{\perp}}
=
(e^{c_1 \delta+c_2 a})^{d_{\alpha}}
\left.\frac{\p^a x^{\alpha}}{\p z^a}
\right|_{V_{c_1,c_2}^{\perp}}
=
\left.\frac{\p^a x^{\alpha}}{\p z^a}
\right|_{V_{0,0}^{\perp}}
=0, 
\end{eqnarray*}
where we used 
$$
\frac{\p}{\p \widetilde{z}^{i}}=
\frac{\p}{\p z^{i}}
$$
for $1 \leq i \leq n$, (\ref{230817.2}) and Proposition \ref{230817.3}\hskip0.2mm(iv). 
By the same argument, we have 
$$
\left.\frac{\p \widetilde{x}^{\alpha}}{\p \widetilde{z}^{\beta}}
\right|_{V_{c_1,c_2}^{\perp}}
=
\left.\frac{\p x^{\alpha}}{\p z^{\beta}}
\right|_{V_{0,0}^{\perp}}
$$
and they are constants. 
Then by Proposition \ref{230817.3}\hskip0.2mm(iv), we have the result. 
\qed
\end{pf}
\begin{prop}\label{240309.2}
If the codimension of an elliptic root system is $1$, 
then the $\CC$-span of good basic invariants 
for the admissible triplet $(\widetilde{c}^{ss},\zeta,L)$ 
does not depend on the choice of 
$L \in \widetilde{X}^{\widetilde{c}^{ss}} \cap \widetilde{X}^{reg}$ 
under the assumption that $L$ is of zero type.
\end{prop}
\begin{pf}
If the codimension of an elliptic root system is $1$, 
$X^{\widetilde{c}^{ss}}$ is one point 
by the results of \cite[p44]{extendedII} and (\ref{240309.4}).
Then $X^{\widetilde{c}^{ss}} \cap X^{reg}$ is also 
one point by Proposition \ref{230402.3}. 
Thus we have the result by Proposition \ref{240309.5}.
\qed
\end{pf}

\section{Good basic invariants for codimension $1$ cases}
In this section, we consider the cases of the elliptic root systems 
of codimension $1$ (i.e. $d_{n-1}<d_n$). 
We fix an admissible triplet $(g,\zeta,L)$ with $L$ of zero type. 
\subsection{Admissible triplet for codimension $1$ cases}
We put
\begin{equation}\label{240308.3}
d_{0}=0.
\end{equation}
By the codimension $1$ assumption that $d_{n-1}<d_n$, we have a duality:
\begin{equation}
d_{\alpha}+d_{n-\alpha}=d_n\quad (0 \leq \alpha \leq n). \label{200314.2}
\end{equation}

For the admissible triplet $(g,\zeta,L)$ with $L$ of zero type, we have 
\begin{equation}\label{240308.2}
F=\mathrm{rad}\,I \oplus (L \cap F).
\end{equation}
By (\ref{240308.2}), we see that $\widetilde{I}$ on $L \cap F$ is nondegenerate. 
Then the orthogonal complement 
\begin{equation}\label{200322.10}
(L \cap F)^{\perp}
:=
\{x \in L\,|\,
\widetilde{I}(x,y)=0\ \forall y \in L \cap F\ \}
\end{equation}
of $L \cap F$ gives a direct decomposition of $L$:
\begin{equation}
L=(L \cap F)\oplus (L \cap F)^{\perp}. \label{200323.6}
\end{equation}
We have the following proposition.
\begin{prop}\label{200323.5}
%By the decomposition 
%$$
%\widetilde{F}=\mathrm{rad}\,I \oplus (L \cap F)\oplus (L \cap F)^{\perp}
%$$
%obtained by $(\ref{200323.7})$ and $(\ref{200323.6})$, 
There exists a basis $z^0,z^1,\cdots,z^n$ 
of $(L \oplus \RR \delta)\otimes_{\RR}\CC$ 
such that 
$z^0=\delta/(-2\pi\sqrt{-1})$, $z^n \in (L \cap F)^{\perp}\otimes_{\RR}\CC$, 
$z^1,\cdots,z^{n-1}\in (L \cap F)\otimes_{\RR}\CC$ with 
\begin{eqnarray}
&&\widetilde{I}(z^{\alpha},z^{\beta})=\delta_{\alpha+\beta,n}
\quad (0 \leq \alpha,\beta \leq n),\label{200310.11}\\
&&g \cdot z^{\alpha}=\zeta^{d_{\alpha}}z^{\alpha}
\quad (0 \leq \alpha \leq n).\label{240307.12}
\end{eqnarray}
\end{prop}

\begin{pf}
Since $g \in O(\widetilde{F},F,\mathrm{rad}\,I)$ and $g$ acts on $L$, 
$g$ acts on $L \cap F$ and on its orthogonal complement 
$(L \cap F)^{\perp}$ defined in (\ref{200322.10}). 
We study the eigenvalues of $g$ on these spaces. 

First we show that the eigenvalue of $g$ on the $1$-dimensional space 
$(L \cap F)^{\perp}$ is $\zeta^{d_n}=1$. 
Since $\mathrm{dim}(L \cap F)^{\perp}=1$, we take 
$0 \neq \xi \in (L \cap F)^{\perp}$. 
Then $g \xi=c_0\xi$ for some $c_0 \in \RR$. 
Since $\xi \in \widetilde{F} \setminus F$, 
$\widetilde{I}(\delta,\xi)\neq 0$. 
By $\widetilde{I}(g\cdot \delta,g \cdot \xi)=\widetilde{I}(\delta,\xi)$, 
we have $c_0=1=\zeta^{d_n}$. 

Then by the condition of admissibility (Definition \ref{200310.1}\hskip0.2mm(iii)), 
the eigenvalues of $g$ on $L \cap F$ are $\zeta^{d_1},\cdots,\zeta^{d_{n-1}}$. 

By the duality (\ref{200314.2}), 
we could take $z^1,\cdots,z^{n-1} \in (L \cap F)\otimes_{\RR}\CC$ such that 
$g \cdot z^{\alpha}=\zeta^{d_{\alpha}}z^{\alpha}$ and 
$\widetilde{I}(z^{\alpha},z^{\beta})
=\delta_{\alpha+\beta,n}$ 
for $1 \leq \alpha,\beta \leq n-1$. 

Take $z^0:=\delta/(-2\pi\sqrt{-1})$. 
Take $z^n \in (L \cap F)^{\perp}\otimes_{\RR}\CC$ such that $\widetilde{I}(z^0,z^n)=1$. 
Since the signature of $L$ is $(n-1,1,0)$, 
we have $\widetilde{I}(z^n,z^n)=0$. 

Then we see that $z^0,\cdots,z^n$ satisfy the conditions 
(\ref{200310.11}) and (\ref{240307.12}). 
\qed\end{pf}

\subsection{Bilinear form and Euler field}
Let $z^0,\cdots,z^n$ be the same as in Proposition \ref{200323.5}. 
Then the set $z^0,\cdots,z^n$ forms a coordinate system 
of $Y$. Then the symmetric $\RR$-bilinear form $\widetilde{I}$ on $\widetilde{F}$ 
defines a $\CC$-bilinear form on $\widetilde{F}\otimes_{\RR}\CC$ and 
this gives
\begin{equation}\label{240308.1}
\widetilde{I}:\Omega(Y)\otimes_{F(Y)}\Omega(Y) \to F(Y)
\end{equation}
by $\widetilde{I}(dz^{\alpha},dz^{\beta})
=\widetilde{I}(z^{\alpha},z^{\beta})$. 

We have 
\begin{equation}\label{240308.4}
\frac{\p}{\p z^n}=\widetilde{I}(dz^0), 
\end{equation}
where $\widetilde{I}:\Omega(Y) \to Der(Y)$ is an isomorphism 
induced by (\ref{240308.1}) 
and $Der(Y)$ is the module of derivations of $F(Y)$. 

The Euler field $E$ defined in (\ref{200322.1}) is interpleted as 
\begin{equation}
E=\frac{(I_R:I)d_n}{m_{max}}\widetilde{I}(d\delta),
\end{equation}
where we see $\delta$ as a function on $Y$. 
Then we have 
\begin{equation}\label{200323.1}
\frac{\p}{\p z^n}=\widetilde{I}(dz_0)=
\frac{1}{-2\pi\sqrt{-1}}
\frac{m_{max}}{(I_R:I)d_n}E
\end{equation}
by $z_0=\frac{\delta}{-2\pi\sqrt{-1}}$. 
\subsection{Bilinear form and Euler field on $S^W$}

Let $x^1,\cdots,x^n$ be a set of basic invariants with 
degrees $d_1 \leq \cdots \leq d_{n-1}\leq d_n$. 
We put 
\begin{equation}\label{200323.3}
x^0:=\delta/(-2\pi\sqrt{-1}).
\end{equation}

The Euler field $E$ satisfies 
$$
E x^{\alpha}=d_{\alpha}x^{\alpha}\quad
(0 \leq \alpha \leq n), 
$$
where $d_0=0$ (see (\ref{240308.3})). 
Then the Euler field $E$ descends to 
\begin{equation}\label{200323.2}
E=\sum_{\alpha=0}^n d_{\alpha}x^{\alpha}\frac{\p}{\p x^{\alpha}}
: \Omega_{S^W} \to S^W, 
\end{equation}
where $\Omega_{S^W}$ is the module of K\"ahler differentials of $S^W$ over $\CC$. 
We define the normalized Euler field $E_{\mathrm{norm}}$ by 
\begin{equation}\label{200322.2}
E_{\mathrm{norm}}:=\frac{1}{d_n}E: \Omega_{S^W} \to S^W. 
\end{equation}
Let $z^0,\cdots,z^n$ be the same as in Proposition \ref{200323.5}, 
which form a coordinate system of $Y$. 

By the $W$-invariance of $\widetilde{I}$, 
we have the $S^W$-bilinear form 
\begin{equation}
\widetilde{I}_W:\Omega_{S^W} \otimes_{S^W}\Omega_{S^W} \to S^W
\end{equation}
defined by 
\begin{equation}
\widetilde{I}_W(dx^{\alpha},dx^{\beta})=
\sum_{\gamma_1,\gamma_2=0}^{n}
\frac{\p x^{\alpha}}{\p z^{\gamma_1}}
\frac{\p x^{\beta}}{\p z^{\gamma_2}}
\widetilde{I}(z^{\gamma_1},z^{\gamma_2}) \in S^W
\end{equation}
for $0 \leq \alpha,\beta \leq n$. 

By (\ref{200323.1}) and $x^0=z^0$, 
$\widetilde{I}_W(dx^0):\Omega_{S^W} \to S^W,\omega \mapsto 
\widetilde{I}_W(dx^0,\omega)$ gives 
\begin{equation}\label{200323.9}
\widetilde{I}_W(dx^0)=
\frac{1}{-2\pi\sqrt{-1}}
\frac{m_{max}}{(I_R:I)d_n}E.
\end{equation}

\subsection{Property of a set of basic invariants for codimension $1$}
\begin{prop}
For a set of basic invariants $x^1,\cdots,x^n$, 
\begin{eqnarray}
&&x^1|\rest=\cdots=x^{n-1}|\rest=0,\label{200322.11}\\
&&x^n(q) \neq 0\quad(\forall q \in L^{\perp}).
\label{200310.8}
\end{eqnarray}
\end{prop}
\begin{pf}
As for (\ref{200322.11}), it is shown by (\ref{200310.2}). 
We show (\ref{200310.8}). 
Let $x^0=\delta/(-2\pi\sqrt{-1})$ defined in (\ref{200323.3}). 
For any $\alpha,\beta\ (0 \leq \alpha,\beta \leq n)$, we put 
$$
a^{\alpha,\beta}:=\widetilde{I}_W(dx^{\alpha},dx^{\beta}).
$$
Then by 
$\widetilde{I}_W(dx^{\alpha},dx^{\beta})
=\sum_{\gamma_1,\gamma_2=0}^{n}
\frac{\p x^{\alpha}}{\p z^{\gamma_1}}
\frac{\p x^{\beta}}{\p z^{\gamma_2}}
\widetilde{I}(z^{\gamma_1},z^{\gamma_2})$, 
we have 
$$
\mathrm{det}(a^{\alpha,\beta})_{0 \leq \alpha,\beta \leq n}
=\mathrm{det}(\frac{\p x^{\alpha}}{\p z^{\gamma_1}})_{0 \leq \alpha,\gamma_1 \leq n}
\mathrm{det}(\frac{\p x^{\beta}}{\p z^{\gamma_2}})_{0 \leq \beta,\gamma_2 \leq n}
\mathrm{det}(\widetilde{I}(z^{\gamma_1},z^{\gamma_2}))_{0 \leq \gamma_1,\gamma_2 \leq n}.
$$
By Proposition \ref{200322.12}\hskip0.2mm(ii), 
$\mathrm{det}(a^{\alpha,\beta})_{0 \leq \alpha,\beta \leq n}$ is not 
$0$ for any $q \in L^{\perp}$. 

On the other hand, 
$\mathrm{det}(a^{\alpha,\beta})_{0 \leq \alpha,\beta \leq n} \in S^W_{d_n(n+1)}$. 
Thus we could expand 
\begin{equation}
\mathrm{det}(a^{\alpha,\beta})_{0 \leq \alpha,\beta \leq n}=\sum_{j=0}^{n+1} A_j(x^n)^j
\label{200403.1}
\end{equation}
for $A_j \in F(H)[x^1,\cdots,x^{n-1}] \cap S^W_{d_n(n+1-j)}$. 
Evaluating the RHS of (\ref{200403.1}) at $q \in L^{\perp}$, we have 
$$
A_{n+1}(q)(x^n(q))^{n+1}
$$
which is not zero. Then we have $x^n(q) \neq 0$. 
\qed\end{pf}
The following proposition is a key observation in the study of 
good basic invariants. 
\begin{prop}\label{200322.9}
For $x^0=\delta/(-2\pi\sqrt{-1})$ defined in $(\ref{200323.3})$
and a set of basic invariants $x^1,\cdots,x^n$, 
the following conditions are equivalent.
\begin{enumerate}
\item $\frac{\p x^n}{\p z^n}|\rest$ is a nonzero constant.
\item $x^n|\rest$ is a nonzero constant.
\item $\widetilde{I}_W(dx^n,dx^n)|\rest=0$.
\item $(\p/\p x^n)^2\widetilde{I}_W(dx^n,dx^n)=0$.
\end{enumerate}
\end{prop}
\begin{pf}
By the equation (\ref{200323.1}), we have 
\begin{equation}\label{240307.11}
\frac{\p x^n}{\p z^n}=
\frac{1}{-2\pi\sqrt{-1}}
\frac{m_{max}}{(I_R:I)d_n}E x^n=
\frac{1}{-2\pi\sqrt{-1}}
\frac{m_{max}}{(I_R:I)}x^n. 
\end{equation}
Then (i) is equivalent to (ii). 

By $\widetilde{I}_W(dx^n,dx^n)=
\sum_{\gamma_1,\gamma_2=0}^{n}
\frac{\p x^n}{\p z^{\gamma_1}}
\frac{\p x^n}{\p z^{\gamma_2}}
\widetilde{I}(z^{\gamma_1},z^{\gamma_2})$ 
, $\left.\frac{\p x^n}{\p z^{\gamma_1}}\right|_{L^{\perp}}=0$ if $\gamma_1 \neq 0,n$ 
and $\widetilde{I}(z^{\gamma_1},z^{\gamma_2})=\delta_{\gamma_1+\gamma_2,n}$, 
we have 
$$
\left.\widetilde{I}_W(dx^n,dx^n)\right|_{L^{\perp}}=
2
\left.\frac{\p x^n}{\p z^{0}}\right|_{L^{\perp}}
\left.\frac{\p x^n}{\p z^{n}}\right|_{L^{\perp}}.
$$
By (\ref{240307.11}) and 
$\left.\frac{\p x^n}{\p z^{0}}\right|_{L^{\perp}}
=\frac{\p (x^n|\rest)}{\p z^{0}}$, we have 
$$
\left.\widetilde{I}_W(dx^n,dx^n)\right|_{L^{\perp}}=
2
\frac{\p (x^n|\rest)}{\p z^{0}}
\left(\left.
\frac{1}{-2\pi\sqrt{-1}}
\frac{m_{max}}{(I_R:I)}x^n
\right|_{L^{\perp}}\right). 
$$
By (\ref{200310.8}), 
the condition (ii) is equivalent to 
the condition (iii). 

Since $\widetilde{I}_W(dx^n,dx^n) \in S^W_{2d_n}$, 
we have 
$$
\widetilde{I}_W(dx^n,dx^n)=A (x^n)^2+B (x^n)+C
$$
with $A \in S^W_0=F(H)$, $B \in F(H)[x^1,\cdots,x^{n-1}]\cap S^W_{d_n}$, 
$C \in F(H)[x^1,\cdots,x^{n-1}]\cap S^W_{2d_n}$. 
By $x^{\alpha}|\rest=0$ for $1 \leq \alpha \leq n-1$, 
we have 
$$
\widetilde{I}_W(dx^n,dx^n)|\rest=A (x^n|\rest)^2. 
$$
By (\ref{200310.8}), (iii) is equivalent to $A=0$ and 
it is equivalent to (iv). 
\qed\end{pf}

\subsection{Good basic invariants and the bilinear form}
\begin{thm}\label{230403.1}
For an admissible triplet $(g,\zeta,L)$ with $L$ of zero type 
for the elliptic root system of codimension $1$, 
let $z^0,\cdots,z^n$ be the basis of 
$(L \oplus \RR \delta)\otimes_{\RR}\CC$ 
defined in Proposition $\ref{200323.5}$. 
Put $x^0=\delta/(-2\pi\sqrt{-1})$ defined in $($\ref{200323.3}$)$. 
Let $x^1,\cdots,x^n$ be a set of good basic invariants 
compatible with the $g$-homogeneous basis $z^1,\cdots,z^n$ of $L\otimes_{\RR}\CC$. 
\begin{enumerate}
\item 
\begin{equation}\label{200323.14}
x^n|\rest=
\frac{(-2\pi\sqrt{-1})(I_R:I)}{m_{max}}.
\end{equation}
\item Any 
$\widetilde{I}_W(dx^{0},dx^{\beta})$ $(\beta=0,\cdots,n)$ 
is written as follows:
\begin{equation}\label{200323.8}
\widetilde{I}_W(dx^{0},dx^{\beta})=
\frac{m_{\beta}}{(-2\pi\sqrt{-1})(I_R:I)}x^{\beta}.
\end{equation}
\item 
Any $\widetilde{I}_W(dx^{\alpha},dx^{\beta})$ $(\alpha,\beta=1,\cdots,n)$ 
is written by Taylor coefficients  
\begin{equation}\label{200318.5}
\left.\frac{\partial^a x^{\alpha}}
{\partial z^a}\right|_{L^{\perp}},\ \,
\left.
\left(
\frac{\p}{\p z^0}\frac{\partial^a x^{\alpha}}
{\partial z^a}
\right)\right|_{L^{\perp}}
\ 
(\one,\,a \in \multi,\, d\cdot a=d_{\alpha}+d_n)
\end{equation}
as follows:
\begin{eqnarray}
&&\widetilde{I}_W(dx^{\alpha},dx^{\beta})\nonumber\\
&=&\delta_{\alpha+\beta,n}
\left(\frac{1}{x^n|_{L^{\perp}}}\right)
x^{n}
+\sum_{\substack{b=(b_1,\cdots,b_n),\\ b_n=0,\\ d\cdot b=d_{\alpha}+d_{\beta}}}
\frac{1}{b!}
\left.
\left[
\frac{\partial^b}{\partial z^b}
\left(
\frac{\partial x^{\alpha}}{\partial z^{\beta*}}
+
\frac{\partial x^{\beta}}{\partial z^{\alpha*}}
\right)
\right] 
\right|_{L^{\perp}}
x^b,\label{200318.6}
\end{eqnarray}
where $\alpha*=n-\alpha \quad (0 \leq \alpha \leq n)$. 
\end{enumerate}
\end{thm}
\begin{pf}
(i) By the restriction of (\ref{240307.11}) to $L^{\perp}$, 
we have (\ref{200323.14}) by $\p x^n/\p z^n|\rest=1$. 

(ii) As for (\ref{200323.8}), it is obtained by (\ref{200323.9}) 
and (\ref{240311.1}). 

(iii) We prove (\ref{200318.6}). 
By (i) and Proposition \ref{200322.9}, 
$(\p/\p x^n)^2\widetilde{I}_W(dx^n,dx^n)=0$. 
Then for any $\alpha,\beta\ (\two)$, 
$\widetilde{I}_W(dx^{\alpha},dx^{\beta}) \in S^W_{d_{\alpha}+d_{\beta}}$ 
is represented as  
\begin{equation}\label{200309.3}
\widetilde{I}_W(dx^{\alpha},dx^{\beta})
=\sum_{\substack{a=(a_1,\cdots,a_n) \in \multi,\\
 a_n=0,\\ d\cdot a=d_{\alpha}+d_{\beta}-d_n}}
A_a^{\alpha,\beta}
x^ax^{n}
+
\sum_{\substack{b=(b_1,\cdots,b_n)\in \multi,\\
 b_n=0,\\ d\cdot b=d_{\alpha}+d_{\beta}}}
B_b^{\alpha,\beta}
x^b
\end{equation}
for $A_a^{\alpha,\beta},B_b^{\alpha,\beta} \in F(H)$. 

By taking higher order derivatives of the both sides of (\ref{200309.3}) 
with respect to $z^1,\cdots,z^n$ and evaluating them at $L^{\perp}$, 
we determine $A_a^{\alpha,\beta},B_b^{\alpha,\beta}$ 
in the following lemmas. Proofs of these lemmas 
are almost the same as Lemma 6.10--Lemma 6.16 in \cite{handai3}. 
Thus we omit them. 

\begin{lem}
For the cases $d_{\alpha}+d_{\beta} \leq d_n$, 
\begin{equation}
A_a^{\alpha,\beta}=0\ \mbox{if $d_{\alpha}+d_{\beta} <d_n$ or $a \neq 0$}.
\end{equation}
\end{lem}
\begin{lem}
We have 
\begin{equation}
(\mathrm{RHS} \hbox{ of }(\ref{200309.3}))|\rest=
\begin{cases}
A_0^{\alpha,\beta}x^{n}|\rest\ &
\mbox{if $d_{\alpha}+d_{\beta} = d_n$},\\
0\ &\mbox{if $d_{\alpha}+d_{\beta} \neq d_n$}.
\end{cases}
\end{equation}
\end{lem}
\begin{lem}
We have 
\begin{equation}
(\mathrm{LHS} \hbox{ of }(\ref{200309.3}))|\rest
=\delta_{\alpha+\beta,n}. 
\end{equation}
\end{lem}
\begin{lem}
For the cases $d_{\alpha}+d_{\beta}>d_n$, 
take any multi-index $c=(c_1,\cdots,c_n)\in \multi$ such that 
$c_n=0,\ d\cdot c=d_{\alpha}+d_{\beta}-d_n$, we have 
\begin{equation}\label{200309.4}
\left.
\left[
\frac{1}{c!}\frac{\partial^c}{\partial z^c}
(\mathrm{RHS} \hbox{ of }(\ref{200309.3}))
\right]
\right|_{L^{\perp}}
=A_c^{\alpha,\beta}x^{n}|\rest. 
\end{equation}
\end{lem}
\begin{lem}
For the cases $d_{\alpha}+d_{\beta}>d_n$, 
take any multi-index $c=(c_1,\cdots,c_n)\in \multi$ such that 
$c_n=0,\ d\cdot c=d_{\alpha}+d_{\beta}-d_n$, we have 
\begin{equation}\label{200309.5}
\left.
\left[
\frac{1}{c!}\frac{\partial^c}{\partial z^c}
(\mathrm{LHS} \hbox{ of }(\ref{200309.3}))
\right]
\right|_{L^{\perp}}=0. 
\end{equation}
\end{lem}
By these lemmas, for any $\alpha,\beta,c\ (\two,\ c \in \multi)$, we obtain 
\begin{eqnarray}\label{200313.2}
A^{\alpha,\beta}_c=
\begin{cases}
\delta_{\alpha+\beta,n}
\left(\frac{1}{x^{n}|\rest}\right)
&\mbox{if $d_{\alpha}+d_{\beta}=d_n$, $c=0$},\\
0 
&\mbox{otherwise}.
\end{cases}
\end{eqnarray}
\begin{lem}
For any multi-index $c \in \multi$ with 
$c_n=0,\ d\cdot c=d_{\alpha}+d_{\beta}$, we have  
\begin{equation}\label{200309.6}
\left.
\left[
\frac{1}{c!}\frac{\partial^{c}}{\partial z^{c}}
(\mathrm{RHS} \hbox{ of }(\ref{200309.3}))
\right]
\right|_{L^{\perp}}=B_{c}^{\alpha,\beta}. 
\end{equation}
\end{lem}
\begin{lem}
For any multi-index $c \in \multi$ with 
$c_n=0,\ d\cdot c=d_{\alpha}+d_{\beta}$, we have 
\begin{equation}\label{200309.7}
\left.
\left[
\frac{1}{c!}\frac{\partial^c}{\partial z^c}
(\mathrm{LHS} \hbox{ of }(\ref{200309.3}))
\right]
\right|_{L^{\perp}}
=\frac{1}{c!}
\left[
\left.
\left(
\frac{\partial^{c}}{\partial z^{c}}
\frac{\partial x^{\alpha}}{\partial z^{n-\beta}}
\right)\right|_{L^{\perp}}
+
\left.
\left(
\frac{\partial^{c}}{\partial z^{c}}
\frac{\partial x^{\beta}}{\partial z^{n-\alpha}}
\right)\right|_{L^{\perp}}
\right]. 
\end{equation}
\end{lem}
By (\ref{200309.6}) and (\ref{200309.7}), we have 
\begin{equation}
B_c^{\alpha,\beta} 
=
\frac{1}{c!}
\left.
\left[
\frac{\partial^c}{\partial z^c}
\left(
\frac{\partial x^{\alpha}}{\partial z^{n-\beta}}
+
\frac{\partial x^{\beta}}{\partial z^{n-\alpha}}
\right)
\right] 
\right|_{L^{\perp}}
. 
\end{equation}
\qed\end{pf}
\begin{rmk}
We could easily check that 
(\ref{200318.6}) is correct for $\alpha=0$ or $\beta=0$ cases and 
they coincide with (\ref{200323.8}). 
\end{rmk}

\section{Frobenius manifold}
In this section, we discuss the relation 
between the Frobenius structure and 
the good basic invariants for the elliptic root systems 
of codimension $1$. 

\subsection{Frobenius structure }\label{210710}
We assume that the codimension (see Section \ref{240314.2}) of an elliptic root system is $1$. 
This means that the degrees $d_1 \leq \cdots \leq d_n$ 
of a set of basic invariants $x^1,\cdots,x^n$ 
satisfy 
\begin{equation}\label{200314.1}
d_{n-1}<d_n.
\end{equation} 
For the module of $\CC$-derivations $\Der$ of $S^W$, 
the grading  
\begin{equation}
\Der=\bigoplus_{m \in \ZZ}\Der_m,\quad
\frac{\p}{\p x^{\alpha}} \in \Der_{-d_{\alpha}}\quad
(0 \leq \alpha \leq n)
\end{equation}
is induced by the grading of 
$S^W=\bigoplus_{m \in \ZZ}S^W_m$ 
and we see that the lowest degree part is an 
$F(H)$-free module of rank $1$, i.e. we have 
\begin{equation}
\Der_{-d_{n}}=F(H)\frac{\p}{\p x^n}. 
\end{equation}

Under the condition (\ref{200314.1}), 
the Frobenius structure on $S^W$ is constructed by Saito 
\cite{extendedII} and Satake \cite{handai} 
(see also \cite{Hertling}). 

\begin{thm}\label{300.001}
$($
Saito \cite{extendedII}, Satake \cite{handai}
$)$ 
We assume the condition $(\ref{200314.1})$. 
\begin{enumerate}
\item 
There exist 
an $S^W$-nondegenerate symmetric bilinear form $($called the metric$)$ 
$J:\Der \otimes_{S^W} \Der \to {S^W}$, 
an $S^W$-symmetric bilinear form $($called the multiplication$)$ 
$\circ:\Der \otimes_{{S^W}}\Der \to \Der$ on $\Der$ 
and 
a field $e:\Omega_{S^W} \to {S^W}$, 
satisfying the following conditions:
\begin{enumerate}
\item the metric is invariant under the multiplication, 
i.e. $J(X\circ Y,Z)=J(X,Y \circ Z)$ 
for any vector fields $X,Y,Z:\Omega_{S^W} \to {S^W}$, 
\item $($potentiality$)$ 
the $(3,1)$-tensor $\nabla \circ$ is symmetric
$($where $\nabla$ is the Levi-Civita connection of the metric$)$, 
i.e. 
$\nabla_X(Y \circ Z)
-Y\circ \nabla_X (Z)
-\nabla_Y(X \circ Z)
+X\circ \nabla_Y (Z)
-[X,Y]\circ Z=0$, 
for any vector fields $X,Y,Z:\Omega_{S^W} \to {S^W}$, 
\item the metric $J$ is flat, 
\item $e$ is a unit field for $\circ$ and it is flat, i.e. $\nabla e=0$,
\item the Euler field $E_{\mathrm{norm}}$ satisfies 
$Lie_{E_{\mathrm{norm}}}(\circ)=1 \cdot \circ$, and 
$Lie_{E_{\mathrm{norm}}}(J)=1 \cdot J$, 
\item the intersection form coincides with the bilinear form $\widetilde{I}_W$:
$J(E_{\mathrm{norm}},J^*(\omega)\circ J^*(\omega'))
=\widetilde{I}_W(\omega,\omega')$ for 1-forms 
$\omega,\omega' \in \Omega_{S^W}$, 
where $J^*:\Omega_{S^W} \to \Der$ is the isomorphism induced by 
the dual metric $J^*$ of $J$. 
\end{enumerate}
\item Put 
\begin{equation}\label{240303.1}
V:=\{\delta \in \Der_{-d_{n}}\,|\,
Lie_e(Lie_e (\widetilde{I}_W))=0\ \}.
\end{equation}
Then $V$ is $1$-dimensional vector space over $\CC$. 
\item 
Let $(J,\circ,e)$ be a Frobenius structure satisfying the conditions in 
$(\mathrm{i})$. 
Then $e \in V \setminus\{0\}$. 
Conversely for any element $\widetilde{e} \in V \setminus \{0\}$, 
there exists uniquely a Frobenius structure 
$(\widetilde{J},\widetilde{\circ},\widetilde{e})$ satisfying 
the conditions in $(\mathrm{i})$. 
The Frobenius structure 
$(\widetilde{J},\widetilde{\circ},\widetilde{e})$ 
is written as 
$(\widetilde{J},\widetilde{\circ},\widetilde{e})
=(c^{-1}J,c^{-1}\circ,ce)$ for some $c \in \CC^{\times}$. 
%
%
%
%Explicitly $\widetilde{e}$ has a form $\widetilde{e}=ce$ 
%for some $c \in \CC^{\times}$ by $(\mathrm{ii})$. 
%Then .%
%
%
%The unit field $e$ in 
%$(\mathrm{i})$ is a nonzero element of $V$. 
%Conversely for any nonzero element $e$ of $V$, 
%there exist uniquely a Frobenius structure $(J,\circ,e)$ in $(\mathrm{i})$ 
%satisfying the conditions 
%$(\mathrm{i})(\mathrm{a})-(\mathrm{i})(\mathrm{f})$. 
\end{enumerate}
\end{thm}
\begin{pf}
As for (ii), see Saito \cite{extendedII} and 
Satake \cite[Proposition 4.2]{handai}. 

We show (iii). 
For the dual metric $J^*$, we have 
$Lie_e(\widetilde{I}_W)=J^*$ (see \cite[Proposition 7.2]{handai3}) and $Lie_e(J^*)=0$ by 
$Lie_e(J)=0$ (see \cite[p146]{Hertling}). Then we have $e \in V$. 
The remaining parts are shown in \cite[Proposition 5.2]{handai}. 
\qed
\end{pf}
The metric $J$ could be constructed from $c\widetilde{I}_W$ and $e$ 
as follows. 
\begin{prop}\label{200323.15}
For 1-forms $\omega,\omega' \in \Omega_{S^W}$, 
we have 
\begin{equation}
J^*(\omega,\omega')=(Lie_{e}(\widetilde{I}_W))(\omega,\omega').
\end{equation}
\end{prop}
A proof of this proposition is the same as \cite[Proposition 7.2]{handai3}, 
so we omit it. 
\subsection{Frobenius structure via flat basic invariants }
We shall interpret the Frobenius structure by a set of basic invariants $x^1,\cdots,x^n$.
\begin{prop}\label{230808.1}
For $x^0=\delta/(-2\pi\sqrt{-1})$ defined in $(\ref{200323.3})$ 
and a set of basic invariants $x^1,\cdots,x^n$, 
the following conditions are equivalent.
\begin{enumerate}
\item ${\p}/{\p x^n} \in V$. 
\item $(\p/\p x^n)^2\widetilde{I}_W(dx^n,dx^n)=0$. 
\end{enumerate}
\end{prop}
\begin{pf}
We first remark that if $\alpha$ or $\beta$ is not $n$, then 
$$
\widetilde{I}_W(dx^{\alpha},dx^{\beta}) \in \bigoplus_{m=0}^{2d_n-1}S^W_{m}.
$$
Thus $(\p/\p x^n)^2\widetilde{I}_W(dx^{\alpha},dx^{\beta})=0$. 
Then (i) is equivalent to $(\p/\p x^n)^2\widetilde{I}_W(dx^{\alpha},dx^{\beta})=0$ 
for all $0 \leq \alpha,\beta \leq n$ and 
they are equivalent to (ii). 
\qed\end{pf}
%\newpage
%
%
Let $\nabla$ be the connection introduced in Theorem \ref{300.001}. 
By Theorem \ref{300.001}\hskip0.2mm(iii), 
the metric $J$ of the Frobenius structure satisfying 
conditions in Theorem \ref{300.001}\hskip0.2mm(i) 
is unique up to a constant factor. 
Then $\nabla$ and the notion of {\it flatness} 
do not depend on the choice of the Frobenius structures 
in Theorem \ref{300.001}. 
\begin{defn}
A set of basic invariants $x^1,\cdots,x^n$ is called flat 
w.r.t. the Frobenius structure if 
\begin{equation}
\nabla dx^{\alpha}=0\quad
(\one).
\end{equation}
Then $x^0,x^1,\cdots,x^n$ with 
$x^0=\delta/(-2\pi\sqrt{-1})$ defined in (\ref{200323.3}) 
form a flat coordinate system for the Frobenius structure 
(Saito \cite{extendedII}, see also Satake \cite{handai}). 
\end{defn}
We give a description of the multiplication and the metric 
w.r.t. the set of flat basic invariants. 
\begin{prop}\label{200318.4}
We assume that a set of basic invariants $x^1,\cdots,x^n$ with degrees 
$d_1 \leq \cdots \leq d_{n-1}<d_n$ satisfies the conditions of 
Proposition $\ref{230808.1}$. 
Put 
$x^0=\delta/(-2\pi\sqrt{-1})$ defined in $(\ref{200323.3})$. 
Then a set of basic invariants $x^1,\cdots,x^n$ 
is flat with respect to the 
Frobenius strucuture in Theorem $\ref{300.001}$ if and only if 
\begin{equation}\label{200318.3}
\eta^{\alpha,\beta}:=e\widetilde{I}_W(dx^{\alpha},dx^{\beta})\quad
(0 \leq \alpha,\beta \leq n)
\end{equation}
are all elements of $\CC$. 
If a set of basic invariants $x^1,\cdots,x^n$ is flat, then 
the metric $J$ is described by 
\begin{equation}\label{200323.16}
\left(
\eta_{\alpha,\beta}
\right)_{0 \leq \alpha,\beta \leq n}
:=
\left(
J
\left(
\p_{\alpha},\p_{\beta}
\right)
\right)
_{0 \leq \alpha,\beta \leq n}
=
\left(
\eta^{\alpha,\beta}
\right)^{-1}_{0 \leq \alpha,\beta \leq n}
\end{equation}
and the structure constants $C_{\alpha,\beta}^{\gamma}$ of the multiplication 
defined by 
\begin{equation}
\p_{\alpha}{\circ} \p_{\beta}
=\sum_{\gamma=0}^{n}
C_{\alpha,\beta}^{\gamma}
\p_{\gamma} \quad
(0 \leq \alpha,\beta \leq n)
\end{equation}
are described by 
\begin{equation}
C_{\alpha,\beta}^{\gamma}=
\sum_{\alpha',\beta',\gamma'=0}^{n}
\eta_{\alpha,\alpha'}
\eta_{\beta,\beta'}
\p^{\gamma}
\left(
\frac{d_n}{d_{\alpha'}+d_{\beta'}}\widetilde{I}_W(dx^{\alpha'},dx^{\beta'})
\right) 
\end{equation}
for $0 \leq \alpha,\beta \leq n$, $\alpha \neq n$, 
where we denote 
\begin{eqnarray}
\p_{\alpha}=\frac{\p}{\p x^{\alpha}},\quad
\p^{\alpha}=\sum_{\alpha'=0}^{n}
\eta^{\alpha,\alpha'}\frac{\p}{\p x^{\alpha'}}\quad (0 \leq \alpha \leq n).
\end{eqnarray}
\end{prop}

%
%
%By (\ref{200310.9}), the first non-zero part of 
%$\varphi[g,\zeta,q](x^{\alpha})$ $(1 \leq \alpha \leq n)$ is 
%\begin{equation}
%\sum_{b \in \multi,\,d\cdot b=d_{\alpha}}
%\frac{\p^b (x^{\alpha}-x^{\alpha}(q))}{\p z^b}(q)z^{b}
%=\psi[g,\zeta,q](x^{\alpha})=z^{\alpha}.
%\end{equation}
%The next non-zero part is 
%\begin{equation}
%\sum_{b \in \multi,\,d\cdot b=d_{\alpha}+d_n}
%\frac{\p^b (x^{\alpha}-x^{\alpha}(q))}{\p z^b}(q)z^b
%=
%\sum_{b \in \multi,\,d\cdot b=d_{\alpha}+d_n}
%\frac{\p^b x^{\alpha}}{\p z^b}(q)z^b
%.
%\end{equation}
%By these coefficients, 

\begin{pf}
By Proposition \ref{200323.15}, the dual metric of the metric 
of the Frobenius strucrure is constructed 
from the unit $e$ and $\widetilde{I}_W$ by (\ref{200318.3}). 

For the construction of the multiplication from $\widetilde{I}_W$, 
we remind the reader of the notion of the Frobenius potential (see Satake \cite{handai2}). 
The Frobenius potential $F$ is defined by the relation 
\begin{equation}
C_{\alpha,\beta}^{\gamma}=\p_{\alpha}\p_{\beta}\p^{\gamma}F
\quad
(0\leq \alpha,\beta,\gamma \leq n)
\end{equation}
with the structure constants 
$C_{\alpha,\beta}^{\gamma}$ of the product 
%\begin{equation}
%\p_{\alpha}\circ \p_{\beta}=\sum_{\gamma=1}^{n}
%C_{\alpha,\beta}^{\gamma}\p_{\gamma}\quad
%(1 \leq \alpha,\beta \leq n)
%\end{equation}
and it is related with $\widetilde{I}_W$ as 
\begin{equation}
\widetilde{I}_W(dx^{\alpha},dx^{\beta})
=
\frac{d_{\alpha}+d_{\beta}}{d_n}
\p^{\alpha}
\p^{\beta}F\quad
(0 \leq \alpha,\beta \leq n). 
\end{equation}
We remark that if $\alpha \neq n$ and $\eta_{\alpha,\alpha'} \neq 0$, 
then $\alpha' \neq 0$ and $d_{\alpha'}\neq 0$. 
Thus for any $\alpha,\beta,\gamma\ (\alpha \neq n \mbox{ and }
0 \leq \alpha,\beta,\gamma \leq n)$, we have 
\begin{eqnarray}
C_{\alpha,\beta}^{\gamma}
&=&
\p_{\alpha}\p_{\beta}\p^{\gamma}F \nonumber\\
&=&
\sum_{\alpha',\beta'=0}^{n}
\eta_{\alpha,\alpha'}
\eta_{\beta,\beta'}
\p^{\gamma}\p^{\alpha'}\p^{\beta'}F\nonumber\\
&=&
\sum_{\alpha',\beta'=0}^{n}
\eta_{\alpha,\alpha'}
\eta_{\beta,\beta'}
\p^{\gamma}
\left(
\frac{d_n}{d_{\alpha'}+d_{\beta'}}
\widetilde{I}_W(dx^{\alpha'},dx^{\beta'})
\right).
\end{eqnarray}
Then we have the results. 
\qed\end{pf}

\subsection{Good basic invariants and Frobenius structure}
\begin{cor}
For an admissible triplet $(g,\zeta,L)$ with $L$ of zero type
 for the elliptic root system of codimension $1$, 
we have the following results. 
\begin{enumerate}
\item Let $x^0,x^1,\cdots,x^n$ be the same as in Theorem $\ref{230403.1}$. 
Then 
\begin{equation}\label{240314.3}
e=
\left(x^n|_{L^{\perp}}\right)
\frac{\p}{\p x^n}=
\frac{(-2\pi\sqrt{-1})(I_R:I)}{m_{max}}\frac{\p}{\p x^n}
\end{equation}
is an element of $V$. 
Let $J$ be the metric and $\circ$ be the multiplication of a 
unique Frobenius structure with the unit 
$e$ $(\ref{240314.3})$ in 
Theorem $\ref{300.001}\hskip0.2mm(\mathrm{iii})$. 
Then the metric $J$ and the structure constants of the multiplication 
$C_{\alpha,\beta}^{\gamma}$ 
$(0 \leq \alpha,\beta,\gamma \leq n)$ are 
\begin{eqnarray}
&&J(\frac{\p}{\p x^{\alpha}},\frac{\p}{\p x^{\beta}})
=\delta_{\alpha+\beta,n} \quad 
(0 \leq \alpha, \beta \leq n),\label{200323.17}\\
&&C_{\alpha,\beta}^{\gamma}=\label{240317.1}
\begin{cases}
&\frac{\p }{\p x^{\gamma*}}
\left(
\frac{d_n}{d_{\alpha*}+d_{\beta*}}\widetilde{I}_W(dx^{\alpha*},dx^{\beta*})
\right) \quad (\alpha \neq n)\\
&\delta_{\beta,\gamma}\quad (\alpha=n)
\end{cases},
\end{eqnarray}
which are all written 
by Taylor coefficients $(\ref{200318.5})$ 
by $(\ref{200318.6})$. 
\item If a set of basic invariants is good w.r.t. an admissible 
triplet $(g,\zeta,L)$ with $L$ of zero type, 
then it is flat w.r.t. the Frobenius structure of Thoerem $\ref{300.001}$. 
\item The space $\mathrm{Specan}(F(H)[L])=\mathrm{Specan}F(H)[z^0,\cdots,z^n]$
 has a metric induced by the dual metric $\widetilde{I}$ $(\ref{240308.1})$. 
The space $\mathrm{Specan}S^W=\mathrm{Specan}F(H)[x^1,\cdots,x^n]$
 has a metric $J$. 
Then $\psi[g,\zeta,L]:S^W \simeq F(H)[L]$ 
gives the isometry w.r.t. these metric structures. 
\end{enumerate}
\end{cor}
\begin{pf} 
We prove (i). 
For $x^0,x^1,\cdots,x^n$ in (i), they satisfy the conditions 
in Proposition \ref{230808.1} by Proposition \ref{200322.9} and 
Theorem \ref{230403.1}\hskip0.2mm(i). 
By Theorem \ref{300.001}\hskip0.2mm(iii), 
we have a unique 
Frobenius structure with the unit (\ref{240314.3}). 

By Theorem \ref{230403.1}, we have 
\begin{equation}
e\widetilde{I}_W(dx^{\alpha},dx^{\beta})
=\delta_{\alpha+\beta,n} \quad 
(0 \leq \alpha, \beta \leq n).
\end{equation}
By Proposition \ref{200318.4}, a set of $x^1,\cdots,x^n$ is flat and we have 
(\ref{200323.17}) and (\ref{240317.1}). 
By (\ref{200318.6}) in Theorem \ref{230403.1}, (\ref{240317.1}) are all written by 
Taylor coefficients (\ref{200318.5}). 
(ii) is a direct consequence of (i). 
(iii) is a direct consequence of (\ref{200310.11}), 
(\ref{200323.17}) and 
$\psi[g,\zeta,L](x^{\alpha})=z^{\alpha}$ 
for $\one$. 
\qed\end{pf}

\appendix
\section{Non-uniqueness of good basic invariants 
for the case when codimension $>1$} 

In this appendix, we show that the $\CC$-span of the good basic invariants 
depends on the choice of admissible triplets of zero type 
for the case of an elliptic root system of 
type $A^{(1,1)}_1$. 

Let $F$ be an $\RR$-vector space 
defined by $F:=\RR \alpha_1 \oplus \RR \delta \oplus \RR a$. 
Let $R:=\{\pm \alpha_1+m\delta+n a\,|\,m,n \in \ZZ\}$. 
Let $I:F \times F \to \RR$ be a positive semi-definite symmetric 
bilinear form with 
$I(\alpha_1,\alpha_1)=2$ and $\mathrm{rad}\,I=\RR \delta \oplus \RR a$. 
Then $R$ is an elliptic root system of type $A^{(1,1)}_1$ 
belonging to $(F,I)$. 
We put $\widetilde{F}:=F \oplus \RR \Lambda_0$ and 
let $\widetilde{I}:\widetilde{F}\times \widetilde{F} \to \RR$ 
be a symmetric $\RR$-bilinear form such that 
$\widetilde{I}|_F=I$, 
$\widetilde{I}(\Lambda_0,\alpha_1)=
\widetilde{I}(\Lambda_0,a)=
\widetilde{I}(\Lambda_0,\Lambda_0)=0$ and 
$\widetilde{I}(\Lambda_0,\delta)=1$. 
$(\widetilde{F},\widetilde{I})$ gives a hyperbolic extension of $(F,I)$. 

Then an elliptic Weyl group $W$, a Coxeter transformation 
and the domains $Y,H$ are defined. 
The semi-simple part of the Coxeter transformation $\widetilde{c}^{ss}$ 
is identity. 

Put 
\begin{eqnarray*}
L_1&:=&\RR(\alpha_1-\frac{1}{2}a)\oplus \RR \Lambda_0,\\
L_2&:=&\RR(\alpha_1-\frac{1}{2}\delta) \oplus \RR(\Lambda_0+\frac{1}{4}\alpha_1
-\frac{1}{8}\delta).
\end{eqnarray*}
Then we could easily check that 
$(\widetilde{c}^{ss},1,L_1)$ and $(\widetilde{c}^{ss},1,L_2)$ are
admissible triplets of zero type. 

\begin{prop}
Let $x^1,\ x^2$ $($resp. $\widetilde{x}_1,\ \widetilde{x}_2$$)$ 
be a set of good basic invariants for the admissible 
triplet $(\widetilde{c}^{ss},1,L_1)$ $($resp. $(\widetilde{c}^{ss},1,L_2)$$)$. 
Then the $\CC$-span of $x^1,\ x^2$ and 
the $\CC$-span of $\widetilde{x}_1,\ \widetilde{x}_2$ do not coincide. 
\end{prop}
\begin{pf}
We put $x^0=\widetilde{x}^0=z^0=\widetilde{z}^0=\delta/(-2\pi\sqrt{-1})$. 
Let $z^1,\,z^2$ (resp. $\widetilde{z}^1,\,\widetilde{z}^2$) be 
a basis of $L_1$ (resp. $L_2$). 

We assume that the $\CC$-span of basic invariants $x^1,x^2$ and 
the $\CC$-span of basic invariants $\widetilde{x}^1,\widetilde{x}^2$ 
coincide. 
Then $\CC$-span of basic invariants $x^0,x^1,x^2$ and 
the $\CC$-span of basic invariants $\widetilde{x}^0,\widetilde{x}^1,\widetilde{x}^2$ 
coincide. This implies 
$$
\det \left(
\frac{\p x^{\alpha}}{\p z^{\beta}}
\right)_{0 \leq \alpha,\beta \leq 2}
=c
\det \left(
\frac{\p \widetilde{x}^{\alpha}}{\p \widetilde{z}^{\beta}}
\right)_{0 \leq \alpha,\beta \leq 2}
$$
for some $c \in \mathbb{C}^{\times}$ because 
$\widetilde{z}^0,\,\widetilde{z}^1,\,\widetilde{z}^2$ could be obtained by 
affine transformation of 
${z}^0,\,{z}^1,\,{z}^2$. 

By the discussion in the proof of Proposition \ref{200322.12} (ii), we have 
$$
\det \left(
\frac{\p x^{\alpha}}{\p z^{\beta}}
\right)_{1 \leq \alpha,\beta \leq n}
=c
\det \left(
\frac{\p \widetilde{x}^{\alpha}}{\p \widetilde{z}^{\beta}}
\right)_{1 \leq \alpha,\beta \leq n}. 
$$
We prepare the Weyl denominator. 
We put $\Lambda_1:=\Lambda_0+\frac{1}{2}\alpha_1$ and 
$\rho:=\Lambda_0+\Lambda_1$. 
We also put 
$$
\Delta^{+}:=\{\alpha_1+n\delta\,(n \geq 0),\ k\delta,\ 
-\alpha_1+k\delta\,(k \geq 1)\}. 
$$
Then the Weyl denominator of an affine Lie algebra 
of type $A^{(1)}_1$ is defined by 
\begin{equation}
\Theta_A:=e^{\rho}\prod_{\alpha \in \Delta^{+}}
(1-e^{-\alpha}). \label{240313.1}
\end{equation}
By \cite[p245]{Chevalley2}, the Jacobian determinant equals 
the Weyl denominator $\Theta_A$ 
up to a multiplication of the unit of 
$F(H)$. Then we have 
\begin{eqnarray}
&&
\det \left(
\frac{\p x^{\alpha}}{\p z^{\beta}}
\right)_{1 \leq \alpha,\beta \leq n}
=f(\tau)\Theta_A, \label{240312.1}
\\
&&
c\det \left(
\frac{\p \widetilde{x}^{\alpha}}{\p \widetilde{z}^{\beta}}
\right)_{1 \leq \alpha,\beta \leq n}
=f(\tau)\Theta_A \label{240312.2}
\end{eqnarray}
for some $f(\tau) \in F(H)^{\times}$. 

Since the restriction of the LHS of 
(\ref{240312.1}) (resp. (\ref{240312.2})) 
to $L_1^{\perp}$ (resp. ${L_2}^{\perp}$) is a constant function, 
the $\CC$-span of $\displaystyle{\Theta_A|_{L_1^{\perp}}}$ 
coincides with the one of 
$\displaystyle{\Theta_A|_{L_2^{\perp}}}$, i.e. 
\begin{equation}
\CC \Theta_A|_{L_1^{\perp}}=\CC \Theta_A|_{{L_2}^{\perp}}. 
\label{240313.2}
\end{equation}

On the other hand, we have an explicit description of 
$\Theta_A|_{L_1^{\perp}}$ (resp. $\Theta_A|_{{L_2}^{\perp}}$) 
by eliminating $\Lambda_0,\,\alpha_1,\,a$ from 
(\ref{240313.1}) 
by the relations $\alpha_1-\frac{1}{2}a=\Lambda_0=0$ 
(resp. $\alpha_1-\frac{1}{2}\delta=\Lambda_0+\frac{1}{4}\alpha_1-\frac{1}{8}\delta=0$)
and $a=-2\pi\sqrt{-1}$. 
Then we have 
$$
\Theta_A|_{L_1^{\perp}}=\exp(\frac{-2\pi\sqrt{-1}}{4})
\prod_{n\geq 0}(1+q^n)\prod_{k\geq 1}(1-q^n)\prod_{k \geq 1}(1+q^n),
$$
$$
\Theta_A|_{L_2^{\perp}}=q^{-\frac{1}{4}}
\prod_{n\geq 0}(1-q^{n+\frac{1}{2}})
\prod_{k\geq 1}(1-q^n)
\prod_{k \geq 1}(1+q^{n-\frac{1}{2}}),
$$
with notation $q=e^{-\delta}$. 
These contradict to (\ref{240313.2}). 
Thus we have the result. 
\qed
\end{pf}

%\end{large}
\end{document}